\documentclass[11pt, reqno]{amsart}
 
\usepackage[cp1251]{inputenc}
\usepackage{amsmath,amsthm,mathrsfs}
\usepackage{amssymb}
\usepackage[english]{babel}
\usepackage{esint}
\usepackage[titletoc]{appendix}

\usepackage{graphicx}
\usepackage{float}
\usepackage{subfigure}

\usepackage{bbm}
\usepackage{amsfonts,amstext,amssymb,verbatim,epsfig}
\usepackage{pgf,tikz}
\usepackage{float}
\usetikzlibrary{arrows}
\usepackage{hyperref}
\usepackage{cite}
\usepackage{epstopdf}
\usepackage{color}
\usepackage{transparent}
\usepackage{subfigure}
\usepackage{pgfplots}
\usepgfplotslibrary{polar}
\usepgflibrary{shapes.geometric}
\usetikzlibrary{calc}

\usepackage{tikz}

\usepackage{xcolor}
\usepackage[normalem]{ulem}

\hypersetup{colorlinks=true,pdfborder=001,  citecolor  = blue, citebordercolor= {magenta}}
\hypersetup{linkcolor=magenta, linkbordercolor = blue}
\hypersetup{colorlinks,urlcolor=blue}
\makeatletter
\let\reftagform@=\tagform@
\def\tagform@#1{\maketag@@@{(\ignorespaces\textcolor{magenta}{#1}\unskip\@@italiccorr)}}
\renewcommand{\eqref}[1]{\textup{\reftagform@{\ref{#1}}}}
\makeatother

\makeatletter

\DeclareUrlCommand\ULurl@@{%
  \def\UrlLeft{\uline\bgroup}%
  \def\UrlRight{\egroup}}
\def\ULurl@#1{\hyper@linkurl{\ULurl@@{#1}}{#1}}
\DeclareRobustCommand*\ULurl{\hyper@normalise\ULurl@}
\makeatother

\def\lessim{\ \lower4pt\hbox{$
		\buildrel{\displaystyle <}\over\sim$}\ }
\def\gessim{\ \lower4pt\hbox{$\buildrel{\displaystyle >}
		\over\sim$}\ }

\newcommand{\nZ}{\mathit{Z}}
\newcommand{\nY}{\mathit{Y}}

\newcommand{\nA}{\mathit{A}}
\newcommand{\nF}{\mathit{F}}

\newcommand{\dF}{\mathcal{F}}
\newcommand{\sA}{\mathcal{A}}

\newcommand{\rr}{\textbf{r}}

\newtheorem{definition}{\bf Definition}
\newtheorem{theorem}{\bf Theorem}
\newtheorem{lemma}[theorem]{\bf Lemma}
\newtheorem{corollary}[theorem]{\bf Corollary}
\newtheorem{example}{\bf Example}

\theoremstyle{remark}
\newtheorem{remark}{Remark}

\newenvironment{Proof of lemma}{\noindent{\bf Proof of Lemma}}{\hfill$\Box$\newline}
\newenvironment{Proof of theorem}{\noindent{\bf Proof of Theorem}}{\hfill{\footnotesize${\square}$}\newline}
\newenvironment{Proof of theorems}{\noindent{\bf Proof of Theorems}}{\hfill$\Box$\newline}
\newenvironment{Proof of proposition}{\noindent{\bf Proof of Proposition}}{\hfill$\Box$\newline}
\newenvironment{Proof of propositions}{\noindent{\bf Proof of Propositions}}{\hfill$\Box$\newline}
\newenvironment{Proof of exercise}{\noindent{\it Proof of Exercise:}}{\hfill$\Box$}
\begin{document}
\title{The spherical mixed $p-$spin glass at zero temperature}

\author{Yuxin Zhou}
\address{Department of Mathematics, Northwestern Universty}
\email{yuxinzhou2023@u.northwestern.edu}

\begin{abstract}
We consider the spherical mixed $p-$spin models and investigate the structure of the Parisi measure at zero temperature. We prove that for the spherical spin models with $n$ components, the Parisi measure at zero temperature is at most $n-$RSB or $n-$FRSB. We also provide a necessary and sufficient computational criterion for the spherical mixed $p-$spin model to be $k-$RSB or $k-$FRSB, $k \leq n.$

 \end{abstract}

\maketitle
\section{Introduction and main results}\label{main}

Fix $N, n>1$. Consider integers $p_1,\cdots,p_n$ with $2<p_1<p_2 < \cdots < p_n$ and  $\lambda_1,\lambda_2,\cdots, \lambda_{n-1} \in [0,1]$. Set $\lambda_n=1-\sum^{n-1}_{i=1} \lambda_i .$  The Hamiltonian of the spherical mixed $p$-spin model is the Gaussian function defined on the $N$-dimensional sphere  $S_{N}:= \big \{ \sigma \in \mathbb{R}^N: \sum^N_{i=1} \sigma^2_i=N \big\}$ by
\begin{eqnarray*}
H_N(\sigma):=\sum^{n}_{j=1} \sum_{1 \leq i_1,\cdots,i_{p_j} \leq N}  \frac{\sqrt{\lambda_j}}{N^{\frac{p_j-1}{2}}}  g_{i_1,\dots,i_{p_j}} \sigma_{i_1} \cdots \sigma_{i_{p_j}} ,
\end{eqnarray*}
where all $(g_{i_1,\cdots, i_{p_j}})$, $1 \leq i_1,\cdots ,i_{p_j} \leq N$, $j=1,\cdots ,n$  are independent, identically distributed standard Gaussian random variables.

The Gaussian field $H_{N}$ is centered with covariance given by,
\begin{eqnarray*}
\mathbb{E} H_N(\sigma^1) H_N(\sigma^2)=N \xi(R_{1,2}),
\end{eqnarray*}
where  $R_{1,2}:=\frac{1}{N}\sum^N_{i=1} \sigma_i^1 \sigma^2_i$
is the normalized inner product between $\sigma^1$ and $\sigma^2$ and 
\begin{eqnarray}\label{eq:psxi}
\xi(x):=\sum^{n-1}_{j=1} \lambda_j \cdot x^{p_j} +\Big( 1- \sum^{n-1}_{j=1} \lambda_j \Big) \cdot x^{p_n}.
\end{eqnarray}

The partition function $Z_{N} = \int_{S_{N}} \exp\big(\beta H_{N}(\sigma)\big) d \sigma$ and the maximum 
\[
M_{N} = \max_{\sigma \in S_{N}} H_{N}(\sigma) 
\] 
have been studied for a long time by both physics and mathematics community. We refer the readers to the books of Mezard-Parisi-Virasoro \cite{mezard1987spin}, Talagrand \cite{TalagrandVolI}, the recent survey \cite{AMSurvey} and the numerous references therein for results in this direction.

Denote by $\mathcal{K}$ a collection of measures $\nu$ on $[0,1]$ with the form 
\begin{eqnarray*}
\nu(ds)=\mathbbm{1}_{[0,1)} \gamma(x)dx +\Delta \delta_{\{1\}}(dx),
\end{eqnarray*}
where $\gamma(x)$ is a nonnegative, right-continuous, and nondecreasing function on $[0,1)$, $\Delta>0,$ and $\delta_{\{1\}}$ is a Dirac measure at 1.  For $\nu \in \mathcal{K},$ define the Crisanti-Sommers functional\cite{CS} by 
\begin{eqnarray*}
\mathcal{Q}(\nu)=\frac{1}{2} \Big( \int^1_0  \xi'(x) \nu(dx) +\int^1_0 \frac{dx}{\nu \big( (x,1] \big)} \Big).
\end{eqnarray*}
 The Parisi formula for the maximum energy (see \cite[Theorem 1]{ArnabChen15}) states that almost surely
\begin{eqnarray*}
\lim_{N \rightarrow \infty} N^{-1} M_N= \inf_{\nu \in \mathcal{K}} \mathcal{Q}(\nu).
\end{eqnarray*}
The Crisanti-Sommers functional $\mathcal{Q}$ is strictly convex on $\mathcal{K}$ and the right-hand side then has a unique minimizer, which is denoted by 
\begin{eqnarray*}
\nu_P(dx)=\gamma_P(x) \mathbbm{1}_{[0,1)}(x)dx + \Delta_P \delta_{\{1\}}(dx).
\end{eqnarray*}
We denote by $\rho_P$ the measure on $[0,1)$ induced by $\gamma_P$, i.e.,
\begin{eqnarray*}
\gamma_P(x) = \rho_P([0,x]), \forall x \in [0,1),
\end{eqnarray*}
and we call $\nu_P$ the Parisi measure at zero temperature. At any positive temperature $\beta$, the Parisi formula is also valid, that is
\begin{equation}\label{PFpt}
\lim_{N\to \infty} \frac{1}{N} \log Z_{N}(\beta) = \inf_{\mu} \mathcal P_{\beta} (\mu),
\end{equation}
where the infimum in \eqref{PFpt} runs over the space of probability measures on $[0,1]$ and $\mathcal P_{\beta}$ is the Parisi functional described in \cite[Theorem 1.1]{Tal06} and \cite{CS}. The minimizer $\mu_{P}(\beta)$ of \eqref{PFpt} is also unique, which is called the Parisi measure at inverse temperature $\beta$.

For the past decades, Parisi measure has attracted a lot of attention in both physics and mathematics. It is known that the structure of the Parisi measure qualitatively describes
the energy landscape of the corresponding Hamiltonian, and it relates to efficiency of a large class of optimization algorithms \cite{el2021optimization, montanari2019optimization, subag2018following}.  Although the role and importance of the order parameter have been unveiled, the structure of the Parisi measure in many
models is yet to be discovered. 
Recently, a full phase diagram for the spherical models with two components has been characterized\cite{AZhou}. However, when the number of components increase, there are much more possibilities of the structure of the Parisi measure at zero temperature, which needs our deeper investigation. 

In this paper, we consider the spherical spin models with arbitrary number of components. We prove that for the spherical model with $n$ components, the support of the Parisi measure at zero temperature contains at most $n+1$ isolated points. The Parisi measure is then at most $n-$RSB or $n-$FRSB.

We also provide a computational criterion to justify $k-$RSB and $k-$FRSB phases with $k \leq n$ for the spherical model with $n$ components. The computational criterion lowers the difficulty of justifying the Parisi measure at zero temperature by replacing the exact equations with some more orderly inequalities, which are much easier to verify. We list a specific spherical spin model with 3 components as an example to show how the phases of the Parisi measure at zero temperature are determined by the computational criterion.

\begin{figure}[H] 
\centering 
\includegraphics[width=9cm,height=4.9cm]{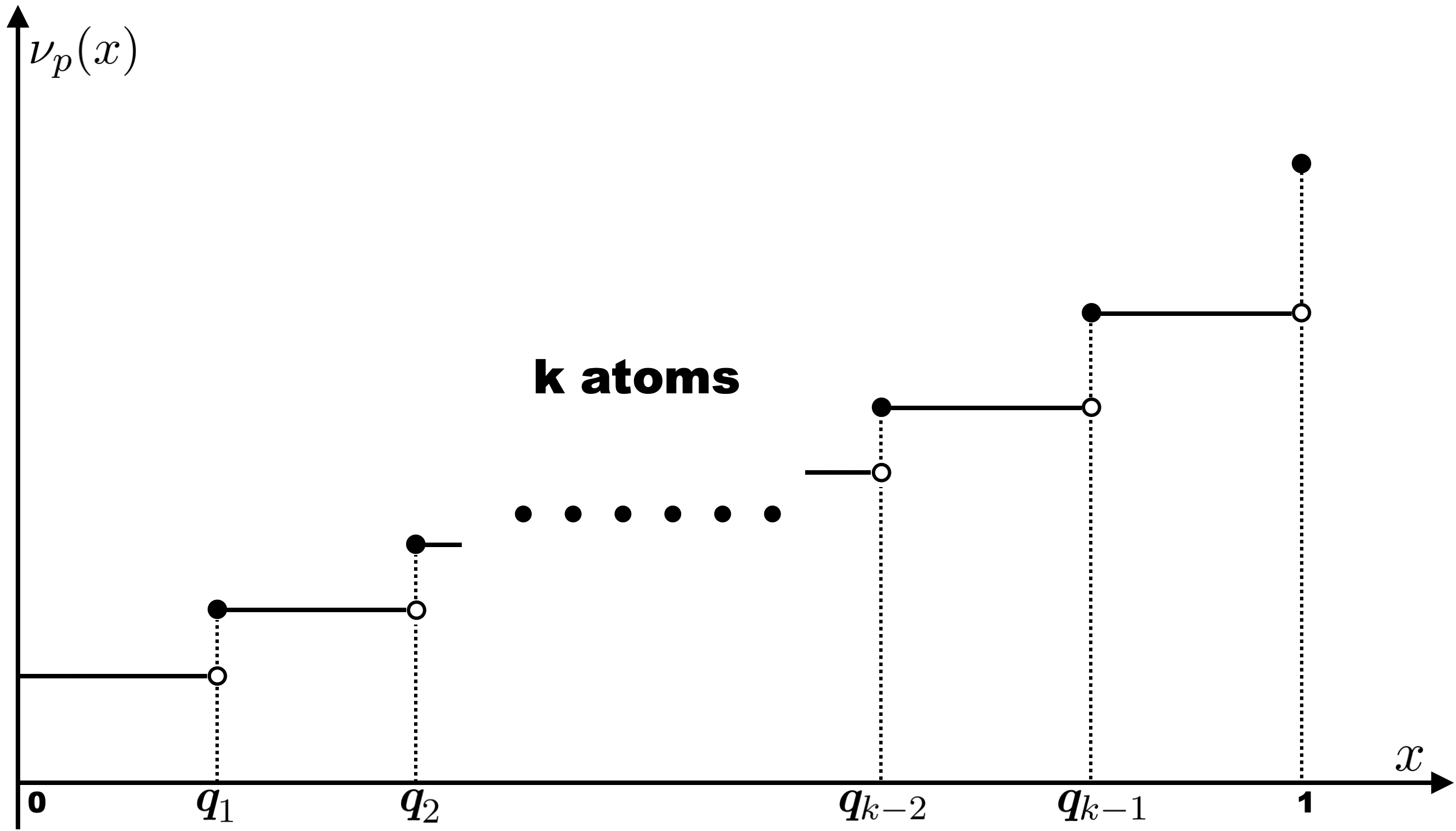}
\caption{The $k$-RSB Parisi measure at zero temperature.} 
\label{kkkrsb}
\end{figure}

We first define the possible structures of the Parisi measure $\nu_P$ at zero temperature for the spherical mixed $p-$spin model as follows. 
\begin{definition}[$k$-RSB] For $k\geq0$, a Parisi measure $\nu_{P}$ at zero temperature has $k$ steps of replica symmetry breaking or, equivalently, the model is called $k$-RSB at zero temperature (Figure \ref{kkkrsb}), if the  support of $\nu_{P}$ is discrete and has exactly $k+1$ points. 
Specifically, there exist nonnegative increasing sequences $\{ m_i \}^k_{i=1}$ and $\{ q_i \}^k_{i=0}$ with $q_0=0,q_k=1$ such that $\nu_P \in \mathcal{K}$ takes the form of 
\begin{eqnarray*}
\nu_P(dx)=\sum^{k}_{i=1} m_i \cdot \mathbbm{1}_{[q_{i-1},q_i)}(x)dx+\Delta \cdot \delta_{\{1\}}(dx).
\end{eqnarray*}
When $k=0$ and then $\nu_P(dx)=0 \cdot \mathbbm{1}_{[0,1)}(x)dx+\Delta \cdot \delta_{\{1\}}(dx)$,  the model is also called Replica Symmetric (RS) at zero temperature.

\end{definition}

\begin{definition}[$k$-FRSB] A Parisi measure $\nu_{P}$ is called  k-FRSB at zero temperature, if for $k>0$, there exist 
\begin{enumerate}
\item $s_1+s_2+\cdots+s_t=k$ for some fixed $t \in [1,k].$ Define $w_j=\sum^j_{i=1} s_i$ for $j=0,1,2,\cdots,t.$ 
\item $q^j_l \in(0,1]$, with $q^{j-1}_{w_{j-1}}<q^j_{w_{j-1}} < q^j_{w_{j-1}+1} <\cdots< q^j_{w_j}$, for $j=1, \cdots, t$ and $l=w_{j-1},w_{j-1}+1,\cdots,w_j$. Set $q^0_0=q^1_0=0$ and $q^{t}_{w_{t}}=1$.
\item $0<m_{1} < m_2 < \cdots < m_k$,
\item a strictly increasing function $w(x)$, 

\end{enumerate}
such that $\nu_P$ takes the form of 
\begin{eqnarray*}
\nu_P(dx)=\sum^{t}_{j=1}\sum^{w_j}_{l=w_{j-1}+1} m_l \cdot \mathbbm{1}_{[q^j_{l-1},q^j_l)}(x)dx+\sum^{t-1}_{j=1} \omega(x) \cdot \mathbbm{1}_{[q^j_{w_j},q^{j+1}_{w_j})}(x)dx+\Delta \cdot \delta_{\{1\}}(dx).
\end{eqnarray*}

\end{definition}

 \begin{remark}
We define the set $\mathcal{F}:= \{ w_j | q^j_{w_j} \neq q^{j+1}_{w_j}, j=1, \cdots t-1 \} $, which marks the location of full replica symmetry breaking for the Parisi measure $\nu_P.$ To be more specific, for $ w_j \in \mathcal{F}$, there is a smooth density between the $w_j-$th and $w_{j+1}-$th atoms, which is the $j-$th smooth density in $\nu_P$. When $t=0$ and  $\mathcal{F}$ is empty, the Parisi measure $\nu_P$ is then $k$-RSB. The following figure \ref{kkkfrsb} gives us an example of the $k$-FRSB Parisi measure at zero temperature.
 \end{remark}
\begin{figure}[H] 
\centering 
\includegraphics[width=11cm,height=5.5cm]{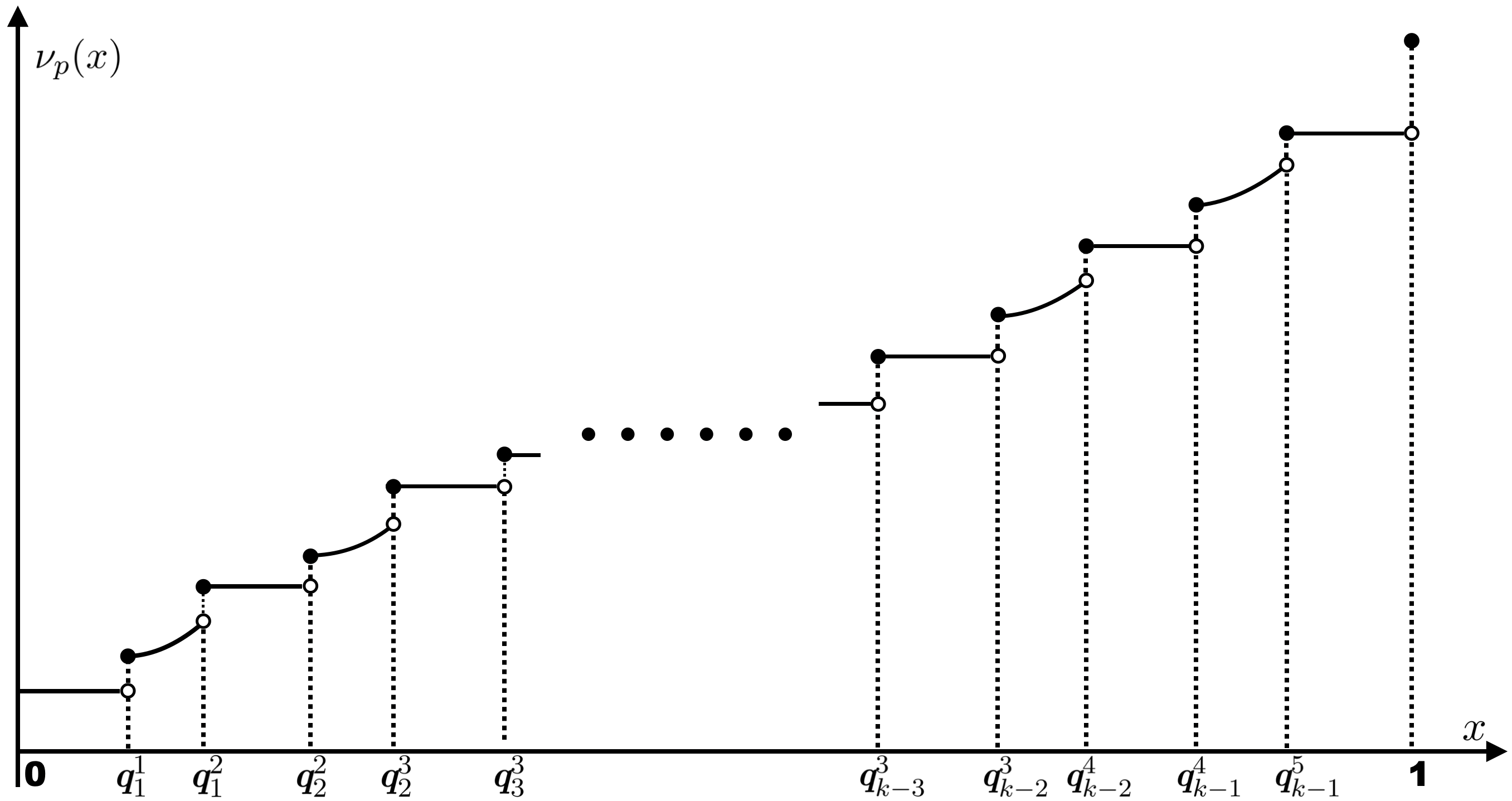}
\caption{An example of the $k$-FRSB Parisi measure at zero temperature with $\mathcal{F}=\{ 1,2,k-2,k-1 \}$. In this case, $s_1=s_2=s_4=s_5=1,s_3=k-4$ and $t=5$. } 
\label{kkkfrsb}
\end{figure}

Based on the two definitions above, the Parisi measure $\nu_P$ at zero temperature is called  at most $n$-RSB or $n$-FRSB if it is $k$-RSB or $k$-FRSB for some $k \leq n.$

Our main results about the Parisi measure $\nu_P$ at zero temperature are as follows:
\begin{theorem}\label{atmostnrsb}
For the spherical model with n components with $\xi(x)=\sum^{n}_{i=1} \lambda_i \cdot x^{p_i}$, the Parisi measure $\nu_P$ at zero temperature is at most $n$-RSB or $n$-FRSB, i.e., the support of $\nu_P$ contains at most $n+1$ isolated points.
\end{theorem}

 For a finite sequence of nonempty subsets $\{ A_i \}_{1 \leq i \leq t }$ with $A_i \subseteq \mathbb{R}^{s_i}$, we define the following relations:
\begin{eqnarray*}
A_1 \tilde{<} A_2 \tilde{<} \cdots \tilde{<} A_{t-1} \tilde{<} A_t ,
\end{eqnarray*}  
if the following conditions are satisfied:
\begin{enumerate}
\item for $i=1,\cdots,t-1$, $\max_{(x_1,\cdots,x_{s_i})\in A_i} x_{s_i} \leq \min_{(x_1,\cdots,x_{s_{i+1}})\in A_{i+1}} x_{1},$
\item $\forall (x_1,\cdots,x_{s_1})\in A_1$, $x_1=0.$ $\forall (x_1,\cdots,x_{s_t})\in A_t$, $x_{s_t}=1.$
\end{enumerate}

We then provide the following criterions to characterize the Parisi measure $\nu_P$ at zero temperature for the model $\xi$ by a sequence of sets $\mathcal{H}$, which will be defined later. We first introduce a criterion for the Parisi measure at zero temperature to be $k-$RSB, for $k \leq n$:
\begin{theorem}\label{kRSB}
The Parisi measure $\nu_P$ for $\xi$ at zero temperature  is $k-$RSB if and only if 
\begin{enumerate}
\item there exists $\bar{x}=(0,x_1,\cdots,x_{k-1},1)$ such that $\bar{x} \in \mathcal{H}^k$.
\item there doesn't exist $\bar{x}=(0,x_1,\cdots,x_{k},1)$ such that $\bar{x} \in \mathcal{H}^{k+1}$.
\end{enumerate}

\end{theorem}

Similarly, we introduce the following criterion for the Parisi measure at zero temperature to be $k-$FRSB, for $k \leq n:$
\begin{theorem}\label{kFRSB}
The Parisi measure $\nu_P$ for $\xi$ at zero temperature is $k-$FRSB if and only if for $1 \leq t \leq k$ and $s_1,s_2,\cdots,s_t \geq 0$ with $s_1+s_2+\cdots+s_t=k,$ 
\begin{enumerate}
\item the relation $\mathcal{H}^{s_1} \tilde{<} \mathcal{H}^{s_2} \tilde{<} \cdots \tilde{<} \mathcal{H}^{s_t} $ holds, 

\item for any $j=1,\cdots, t$, the relations $\mathcal{H}^{s_1} \cdots \tilde{<} \mathcal{H}^{s_{j-1}} \tilde{<}  \mathcal{H}^{s_{j}+1} \tilde{<}  \mathcal{H}^{s_{j+1}} \tilde{<} \cdots \tilde{<} \mathcal{H}^{s_t} $ or $\mathcal{H}^{s_1} \cdots \tilde{<} \mathcal{H}^{s_{j-1}} \tilde{<}  \mathcal{H}^{1} \tilde{<}  \mathcal{H}^{s_{j}} \tilde{<} \cdots \tilde{<} \mathcal{H}^{s_t} $ don't hold. 




\end{enumerate}

\end{theorem}

Indeed when $k=n,$ condition (2) in Theorem \ref{kRSB} is satisfied automatically. Therefore, we obtain the following result for the Parisi measure at zero temperature to be $n$-RSB:
\begin{corollary}\label{nRSB}
The Parisi measure $\nu_P$ for $\xi$ at zero temperature is $n$-RSB if and only if 
 there exists $\bar{x}=(0,x_1,\cdots,x_{n-1},1)$ such that $\bar{x} \in \mathcal{H}^n$.

\end{corollary}

\begin{remark}
In order to show that the Parisi measure is n-RSB, we now just need to find some $\bar{x} \in \mathcal{H}^n$ rather than solving $n$ exact equations, which lowers the difficulty of verification a lot.
\end{remark}

Moreover, when $k=n,$ condition (2) in Theorem \ref{kFRSB} is also satisfied automatically. Therefore, we obtain the following result for the Parisi measure at zero temperature to be $n-$FRSB:
\begin{corollary}\label{nFRSB}
The Parisi measure $\nu_P$ for $\xi$ at zero temperature is $n-$FRSB if and only if for $1 \leq t \leq n$ and $s_1,s_2,\cdots,s_t \geq 0$ with $s_1+s_2+\cdots+s_t=n,$ it holds that $\mathcal{H}^{s_1} \tilde{<} \mathcal{H}^{s_2} \tilde{<} \cdots \tilde{<} \mathcal{H}^{s_t} $.

\end{corollary}

Now we fix $0<k \leq n$ and introduce the definition of the set $\mathcal{H}$ as follows.
We first define the following functions
 \begin{eqnarray*}
h(x,y,z)=\xi(y)-\xi(x)-\xi'(x)(y-x)+\big[\xi'(y)-\xi'(x)\big]  (y-x) \Big( \frac{1}{( z-1 ) } -\frac{z}{(z-1)^2 }  \log ( z ) \Big),
\end{eqnarray*}
\begin{eqnarray*}
r_1(x,y,z)=\frac{[\xi'(z)-\xi'(y)](z-x)}{(z-y)[\xi'(z)-\xi'(x)]} \text{ and } r_2(x,y)=\frac{\xi''(y)(y-x)}{[\xi'(y)-\xi'(x)]}.
\end{eqnarray*}

For any nonnegative integers $i$, we define 
$\iota(i) =  \left\{
\begin{array}{lcl}
&0,& \text{ for } i \text{  even, }  \\
&1,& \text{ for } i \text{  odd. } 
\end{array} \right. $
Also, for $l>0$ and $0\leq x_{l-1} \leq x_l \leq x_{l+1} \leq 1 $, define $\rr_l(x)= r_1(x_{l-1},x_{l+1},x_l)$.

For any $s \geq 0$ and $\bar{x}=(x_0,x_{1},\cdots,x_{s}) \in [0,1]^{s+1},$
we define, 
\begin{eqnarray}\label{dF}
\dF^s_{l}(\bar{x})=r_2\big(x_{l-(-1)^{s-l}},x_{l} \big) \cdot \prod^{\frac{s-l-\iota(s-l)}{2}}_{i=1} \Big(\frac{\rr_{s+1-2i}(x)}{\rr_{s-2i}(x)} \Big)^{(-1)^{s-l}}, \text{ for } 0 \leq l \leq s,
\end{eqnarray}
and 
\begin{eqnarray}\label{sA}
\sA^s_{l}(\bar{x})=r_1\big(x_{l},x_{l-(-1)^{s-l}},x_{l+(-1)^{s-l}}\big) \cdot \prod^{\frac{s-l-\iota(s-l)}{2}}_{i=1} \Big(\frac{\rr_{s-2i}(x)}{\rr_{s+1-2i}(x)} \Big)^{(-1)^{s-l}},\text{ for } 1 \leq l \leq s-1.
\end{eqnarray}

We then define $$\mathcal{Z}^s(x_0,x_1,\cdots,x_s)=\max \Big\{ \max_{0\leq l \leq s,\iota(s-l)=1} \dF^s_{l}(x), \max_{1\leq l \leq s-1,\iota(s-l)=0} \sA^s_{l}(x) \Big \},$$ and 
$$\mathcal{Y}^s(x_0,x_1,\cdots,x_s)=\max \Big\{ \max_{0\leq l \leq s,\iota(s-l)=0} \dF^s_{l}(x), \max_{1\leq l \leq s-1,\iota(s-l)=1} \sA^s_{l}(x) \Big \}.$$
Without confusion, we abbreviate $\mathcal{Z}^s(x_0,x_1,\cdots,x_s)$ and $\mathcal{Y}^s(x_0,x_1,\cdots,x_s)$ as $\mathcal{Z}^s(\bar{x})$ and $\mathcal{Y}^s(\bar{x})$ respectively. 

For any $s=1,\cdots,k$, we define the following condition for $(x_0,x_{1},\cdots,x_{s}) \in [0,1]^{s+1}$:

\textbf{Condition} $\varkappa(s).$
\begin{enumerate}
\item $x_0<x_1<\cdots<x_s,$
\item $ \text{ For } 1 \leq l \leq s,$ 
$\left \{ \begin{array}{lcl}
(-1)^{s-l} h\Big(x_{l-1},x_{l},  {\mathcal{Z}^s(\bar{x})}^{(-1)^{s-l+1}} {\mathcal{F}^s_{l}(\bar{x})}^{-1} \cdot r_{2}(x_{l-1},x_{l})   \Big) \geq 0, \\
(-1)^{s-l} h \Big(x_{l-1},x_{l},  {\mathcal{Y}^s(\bar{x})}^{(-1)^{s-l}} {\mathcal{F}^s_{l}(\bar{x})}^{-1} \cdot r_{2}(x_{l-1},x_{l}) \Big) \leq 0.
\end{array} \right.$
\end{enumerate}

Therefore, for $s=1,\cdots,k$, we define $\mathcal{H}$ as follows:
\begin{eqnarray*} \mathcal{H}^s =\Big \{ \bar{x}=(x_0,x_{1},\cdots,x_{s})\in[0,1]^{s+1} \Big| \text{ } \bar{ x} \text{ satisfies Condition } \varkappa(s) \Big \}.
\end{eqnarray*}





Now we provide several examples to show how the phases of the Parisi measure at zero temperature are determined by the computational criterion introduced above.

\begin{example}[The spherical model with 2 components \cite{AZhou}] 
 Fix $n=2.$ 
 
 Following the notations in \cite{AZhou}, $(x_{1\rightarrow 2},\lambda_{1 \rightarrow 2})$ is the unique solution of  $\left \{ \begin{array}{lcl}
 h\big(0,x,r_1(0,x,1)\big)= 0,\\
h\big(x,1,r_1(x,1,0)\big) = 0,
 \end{array} \right.$
 where $\lambda_{1 \rightarrow 2}$ is the boundary of 1-RSB phase and 2-RSB phase.
 Also, $(x_{2\rightarrow 2F},\lambda_{2 \rightarrow 2F})$ is the unique solution of  $\left \{ \begin{array}{lcl}
 h\big(0,x,r_2(0,x) \big) = 0,\\
  h\big(x,1,r_2(1,x)^{-1} \big) = 0,
 \end{array} \right.$
 where $\lambda_{2 \rightarrow 2F}$ is the boundary of 2-RSB and 2-FRSB phase.
 
  Set $s=2.$ For $0<x<1,$ $\mathcal{Z}^s(x)=\frac{\xi''(x)(1-x)}{[\xi'(1)-\xi'(x)]}$ and $\mathcal{Y}^s(x)=\frac{\xi'(1)x}{\xi'(x)}$.
For $x \in \mathcal{H}^2$, it holds that
\begin{eqnarray*}
\left \{ \begin{array}{lcl}
 h\big(0,x,r_1(0,x,1)\big)\geq 0,h\big(x,1,r_1(x,1,0)\big) \leq 0,\\
 h\big(0,x,r_2(0,x) \big) \leq 0, h\big(x,1,r_2(1,x)^{-1} \big) \geq 0.
 \end{array} \right.
 \end{eqnarray*}
 For $\lambda \in [\lambda_{1 \rightarrow 2},\lambda_{2 \rightarrow 2F}],$ it holds that $0\tilde{\leq} \mathcal{H}^2 \tilde{\leq} 1$ and by Corollary \ref{nRSB}, the Parisi measure at zero temperature is $2$-RSB.
 
Now set $t=2$ and $s_1=s_2=1.$ For $(0,x) \in \mathcal{H}^1,$
it holds that
$$ h\big(0,x,r_2(0,x) \big) \leq 0, $$
Also, for $(x,1) \in \mathcal{H}^1$, it holds that $$h\big(x,1,r_2(1,x)^{-1} \big) \geq 0.$$ If $0 \tilde{\leq}\mathcal{H}^1 \tilde{\leq} \mathcal{H}^1 \tilde{\leq} 1$, by definition, $\max_{(0,x)\in \mathcal{H}^1} x \leq \min_{(x,1)\in \mathcal{H}^1} x. $ For $x_1=\max_{(0,x)\in \mathcal{H}^1} x,$ it holds that $h\big(0,x,r_2(0,x) \big) = 0.$ Also, for $x_2=\min_{(x,1)\in \mathcal{H}^1} x,$ it holds that $h\big(1,x_2,r_2(1,x_2)^{-1} \big) = 0.$  

 For $\lambda \in [\lambda_{2 \rightarrow 2F},\lambda_{2 \rightarrow 1F}],$ it holds that  $0 \tilde{\leq} \mathcal{H}^1 \tilde{<} \mathcal{H}^1 \tilde{\leq} 1$ and supp $\nu_P=\{0,1\} \cup [x_1,x_2].$ By Corollary \ref{nFRSB}, the Parisi measure at zero temperature is then 2-FRSB.  
  For $\lambda \in [\lambda_{2 \rightarrow 1F},\lambda_{2 \rightarrow 1}],$ it holds that $0 \tilde{\leq} \mathcal{H}^1  \tilde{<} \{1\}.$ In this case, supp $\nu_P=\{0\} \cup [x_1,1].$ The Parisi measure at zero temperature is then 1-FRSB.

\end{example}

Now we consider a specific spherical spin model with 3 components. Set $n=k=3$ and $p_1=4,p_2=28,p_3=84.$ With various choices of parameters $\lambda_1$ and $\lambda_2,$ the Parisi measures at zero temperature in the following examples  are in different phases.

\begin{figure}[H] 
\centering 
\includegraphics[width=9cm,height=4cm]{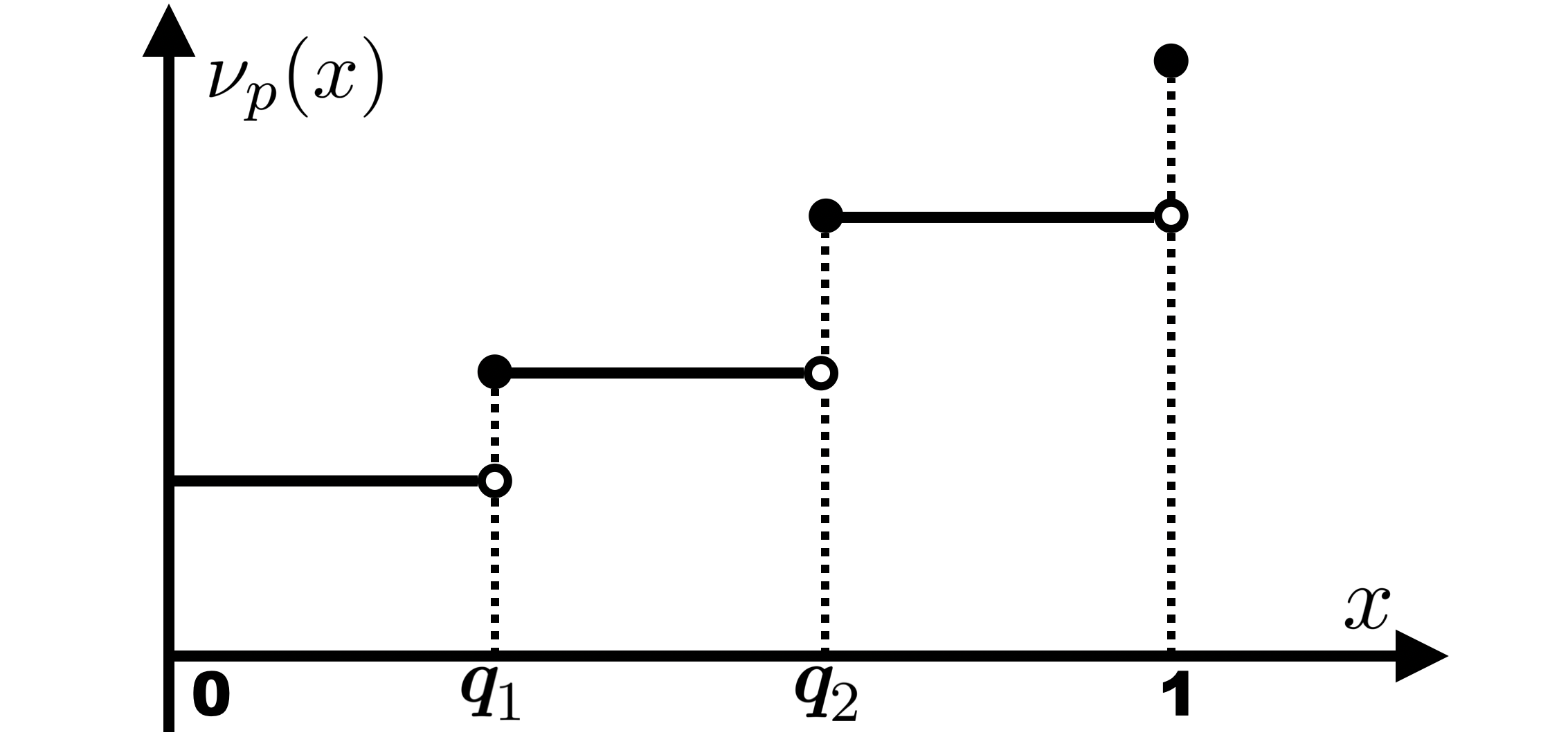}
\caption{The 3-RSB Parisi measure at zero temperature in Example \ref{33rsb}.} 
\label{3rsb}
\end{figure}

\begin{example}[3-RSB]\label{33rsb}
 Set $\lambda_1=0.88$ and $\lambda_2=0.1118.$ Choose $x_1=0.9345$ and $x_2=0.975.$ Then for $(0,x_1,x_2,1) \in [0,1]^4$ satisfies that $\mathcal{Z}^3(x)=\sA_2(x)$ and $\mathcal{Y}^3(x)=\dF_2(x)$. 

It's easy to verify that for $(0,x_1,x_2,1) \in \mathcal{H}^3$,
\begin{eqnarray*}
&&h\Big(0,x_1,\frac{\rr_1(x)}{\rr_2(x)} \cdot \mathcal{Z}^3(x)^{-1}  \Big) > 0 \text{ and } h\Big(0,x_1, \frac{\rr_1(x)}{\rr_2(x)} \cdot \mathcal{Y}^3(x)\Big) < 0, \\
&&h\Big(x_1,x_2, \rr_2(x) \cdot \mathcal{Z}^3(x) \Big) <0 \text{ and } h\big(x_1,x_2,\rr_2(x) \cdot \mathcal{Y}^3(x)^{-1}\big)>0  ,\\
&&h\big(x_2,1, \mathcal{Z}^3(x)^{-1} \big)>0 \text{ and } h\big(x_2,1,\mathcal{Y}^3(x) \big)<0.
\end{eqnarray*}
which yields that $(0,x_1,x_2,1)$ satisfies Condition $\varkappa(3)$ and $(0,x_1,x_2,1) \in \mathcal{H}^3 $. Therefore the Parisi measure at zero temperature is $3$-RSB (Figure \ref{3rsb}) by Corollary \ref{nRSB}.
\end{example}

\begin{figure}[H]
\centering
	\begin{subfigure}
	\centering
	\includegraphics[width=6.5cm,height=4cm]{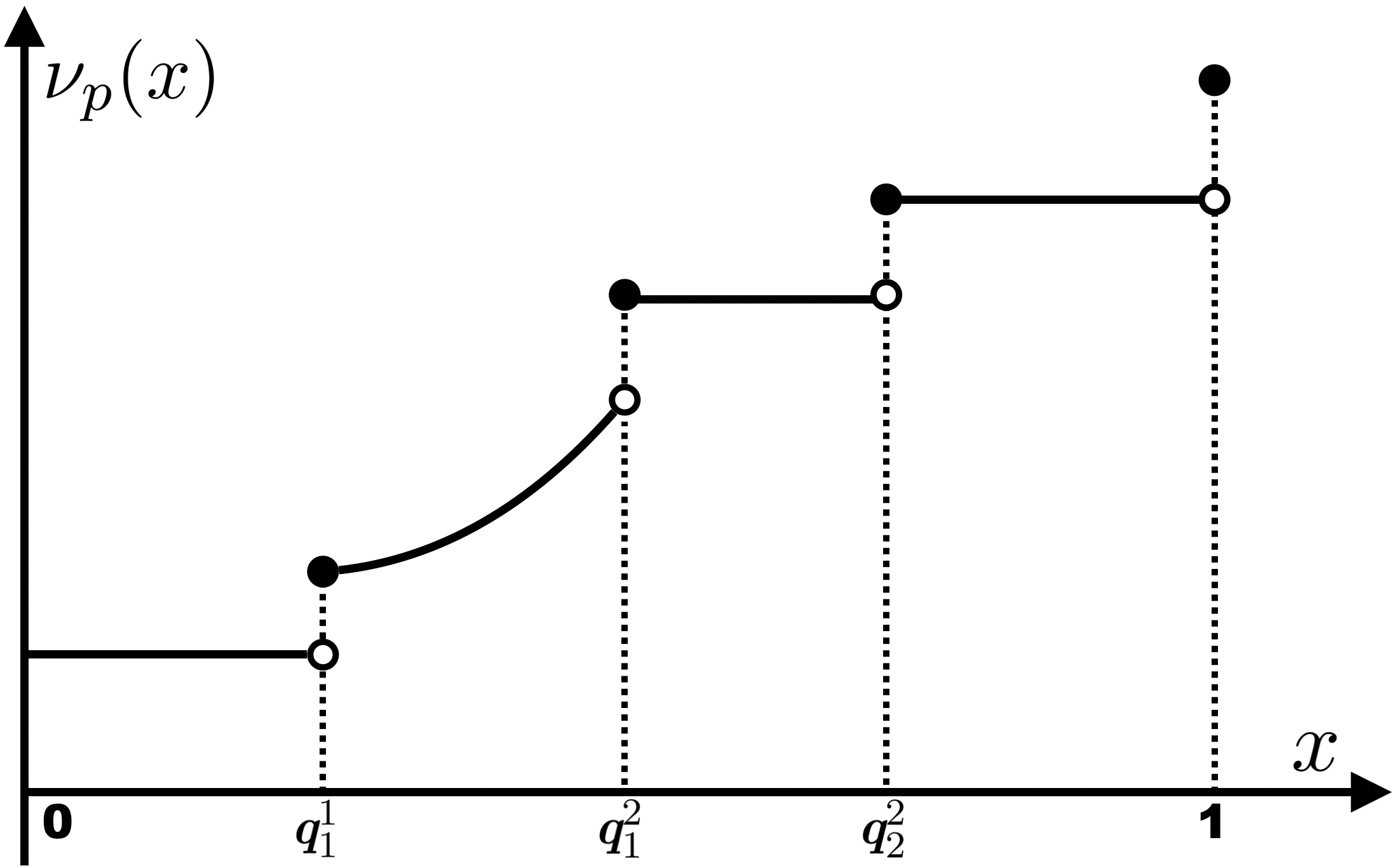} 
\end{subfigure}
\begin{subfigure}
\centering
	\includegraphics[width=6.5cm,height=4cm]{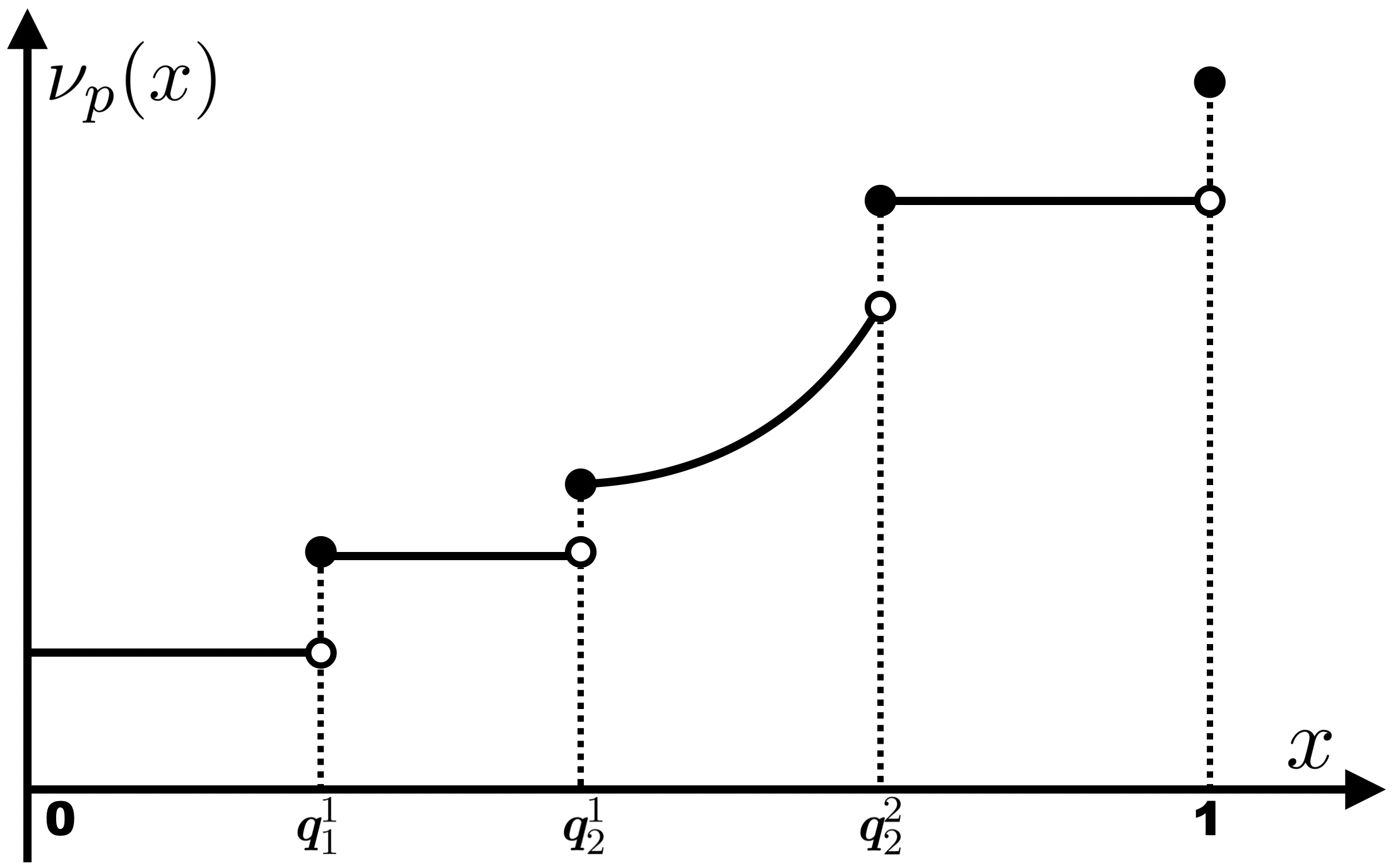} 
	\end{subfigure}
\caption{The Parisi measures at zero temperature in Example \ref{3f1},\ref{3f2}. The subfigure on the left describes the 3-FRSB Parisi measure at zero temperature with $\mathcal{F}=\{ 1\}$. While the subfigure on the right describes the 3-FRSB Parisi measure at zero temperature with $\mathcal{F}=\{ 2\}.$  }
\label{3frsb1+2}
	\end {figure}

Now we denote the following functions with respect to $0<x_1<x_2<1:$
\begin{eqnarray*}
&& \bar{h}^L_1(x_1) = h\big(0,x_1,r_2(0,x_1)\big) \text{ and } \bar{h}^U_1(x_1,x_2) = h\big(0,x_1,\rr_1(x) r_2(x_1,x_2)^{-1}\big),\\
&& \bar{h}^L_2(x_1,x_2) = h\big(x_1,x_2,r_2(x_2,x_1)^{-1} \big) \text{ and } \bar{h}^U_2(x_1,x_2) = h\big(x_1,x_2,r_2(x_1,x_2) \big), \\
&& \bar{h}^L_3(x_1,x_2) = h\big(x_2,1,r_2(x_2,x_1) \rr_2(x) \big) \text{ and } \bar{h}^U_3(x_2) = h\big(x_2,1,r_2(1,x_2)^{-1} \big) .
\end{eqnarray*}
\begin{example}[3-FRSB with $\mathcal{F}=\{ 1 \}$]\label{3f1}
 Set $\lambda_1=0.86$ and $\lambda_2=0.1253.$ 
\begin{enumerate}
\item Set $x_1=0.929$, then $\bar{h}^L_1(x_1)\leq 0$, where $(0,x_1)$ satisfies Condition $\varkappa(1).$ Also set $x'_1=0.93$, then $\bar{h}^L_1(x'_1) \geq 0$, where $(0,x'_1)$ doesn't satisfy Condition $\varkappa(1).$ Moreover, by intermediate value theorem, there exists $q^1_1 \in [0.92,0.93]$ such that $\bar{h}^L_1(q^1_1) = 0$. 
\item Set $x^1_0=0.9352,x^1_1=0.9497$ and $x^2_0=0.936,x^2_1=0.94975.$  Then for $(x_0,x_1) \in \big[x^1_0,x^2_0 \big] \times \big[x^1_1,x^2_1\big],$ it holds that $\mathcal{Z}^2(x)=\dF_1(x)$, $\mathcal{Y}^2(x)=\dF_0(x)$ and 
\begin{eqnarray*}
\bar{h}^U_2(x_0,x_1)  \leq 0, \bar{h}^L_3(x_2) \geq 0.
\end{eqnarray*}
Moreover, for $i=1,2,$
\begin{eqnarray*}
\bar{h}^L_2(x^1_0,x^i_1) \leq 0 , \bar{h}^L_2(x^2_0, x^i_1)  \geq 0,
\end{eqnarray*}
and
\begin{eqnarray*}
\bar{h}^L_3(x^i_0,x^1_1)  \leq 0 ,  \bar{h}^L_3(x^i_0,x^2_1)  \geq 0.
\end{eqnarray*}
Then by intermediate value theorem, there exists $(q^2_1,q^2_2) \in \big[x^1_0,x^2_0 \big] \times \big[x^1_1,x^2_1\big]$ such that  $\left \{ \begin{array}{lcl}
 \bar{h}^L_2 \big(q^2_1,q^2_2 \big)= 0,\\
\bar{h}^L_3(q^2_1,q^2_2) = 0.
 \end{array} \right.$
  Then $ q^2_1= \min_{(x_0,x_1,1)\in \mathcal{H}^2} x_{0}.$
\end{enumerate}

Since $q^1_1 \leq 0.93<0.9352\leq q^2_1,$ then the relation $\mathcal{H}^1 \tilde{<} \mathcal{H}^2$ holds. Then by Corollary \ref{nFRSB}, the Parisi measure is 3-FRSB with $\mathcal{F}= \{ 1  \}$(Figure \ref{3frsb1+2}). In this case, supp $\nu_P=\big \{0,q^2_2,1\big \} \cup \big[q^1_1,q^2_1\big].$
\end{example}

\begin{figure}[H] 
\centering 
\includegraphics[width=9cm,height=5.3cm]{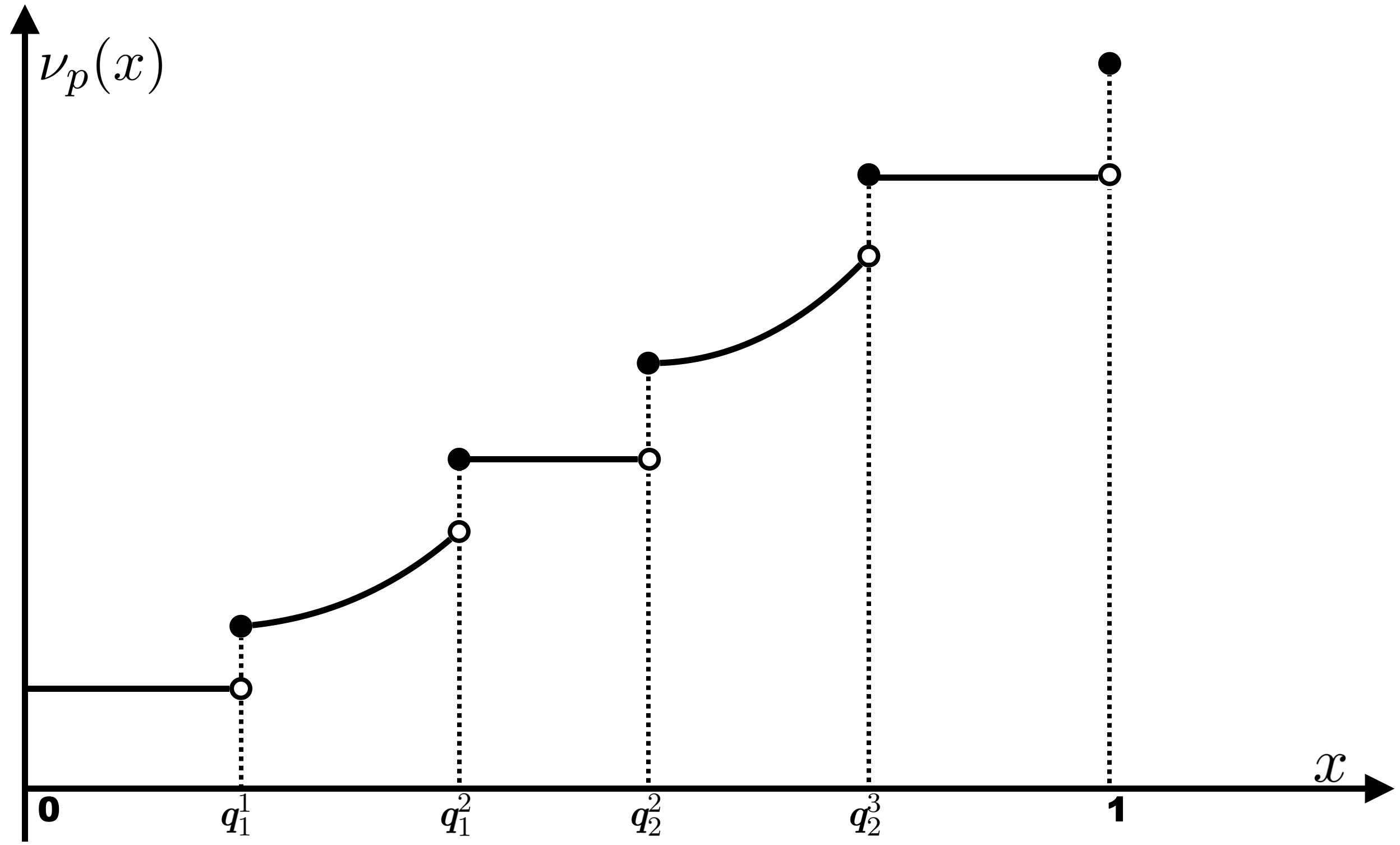}
\caption{The 3-FRSB Parisi measure at zero temperature with $\mathcal{F}=\{ 1,2 \}$ in Example \ref{3f12}. In this case, $t=3$ and $s_1=s_2=s_3=1.$} 
\label{3frsb12}
\end{figure}

\begin{example}[3-FRSB with $\mathcal{F}=\{ 2 \}$]\label{3f2}
Set $\lambda_1=0.88$ and $\lambda=0.1113.$
\begin{enumerate}
\item Set $x^1_1=0.934676,x^1_2=0.97136$ and $x^2_1=0.9346764,x^2_2=0.97139.$  Then for $(x_1,x_2) \in \big[x^1_1,x^2_1 \big] \times \big[x^1_2,x^2_2\big],$ it holds that $\mathcal{Z}^2(x)=\dF_1(x)$, $\mathcal{Y}^2(x)=\dF_2(x)$ and 
\begin{eqnarray*}
\bar{h}^L_1(x_1) \leq 0, \bar{h}^L_2(x_1,x_2)  \geq 0.
\end{eqnarray*}
Moreover, for $i=1,2,$
\begin{eqnarray*}
\bar{h}^U_2(x^i_1,x^1_2 ) \leq 0 , \bar{h}^U_2(x^i_1,x^2_2) \geq 0,
\end{eqnarray*}
and 
$$\left \{ \begin{array}{lcl}
 \bar{h}^U_1(x^1_1,x^1_2) \geq 0 ,  \bar{h}^U_1(x^1_1,x^2_2) \leq 0,\\
 \bar{h}^U_1(x^2_1,x^1_2) \geq 0,  \bar{h}^U_1(x^2_1,x^2_2) \geq 0.
 \end{array} \right.$$
Then by intermediate value theorem, there exists $(q^1_1,q^1_2) \in \big[x^1_1,x^2_1 \big] \times \big[x^1_2,x^2_2\big]$ such that  $\left \{ \begin{array}{lcl}
 \bar{h}^U_1(q^1_1,q^1_2)  = 0,\\
\bar{h}^U_2(q^1_1,q^1_2)  = 0.
 \end{array} \right.$
  Then $ q^1_2= \max_{(0,x_1,x_2)\in \mathcal{H}^2} x_{2}.$

\item Set $x_0=0.9714$, then $\bar{h}^U_3(x_0)  \leq 0$, where $(x_0,1)$ doesn't satisfy Condition $\varkappa(1).$ Also set $x'_0=0.9715$, then $\bar{h}^U_3(x'_0) \geq 0$, where $(x'_0,1)$ satisfy Condition $\varkappa(1).$ Moreover, by intermediate value theorem, there exists $q^2_2 \in [0.9714,0.9715]$ such that $\bar{h}^U_3(q^2_2)  = 0$. 
\end{enumerate}

Since $q^1_2 \leq 0.97139<0.9714\leq q^2_2,$ then the relation $ \mathcal{H}^2 \tilde{<} \mathcal{H}^1$ holds. Then by Corollary \ref{nFRSB}, the Parisi measure is 3-FRSB with $\mathcal{F}= \{ 2  \}$(Figure \ref{3frsb1+2}). In this case, supp $\nu_P=\big \{0,q^1_1,1\big \} \cup \big[q^1_2,q^2_2 \big].$
\end{example}
\begin{example}[3-FRSB with $\mathcal{F}=\{ 1,2 \}$]\label{3f12}
 Set $\lambda_1=0.88$ and $\lambda_2=0.1108.$
\begin{enumerate}
\item Set $x_1=0.9345$, then $ \bar{h}^L_1(x_1)  \leq 0$, where $(0,x_1)$ satisfies Condition $\varkappa(1).$ Also set $x'_1=0.935$, then $ \bar{h}^L_1(x'_1) \geq 0$, where $(0,x'_1)$ doesn't satisfy Condition $\varkappa(1).$ Moreover, by intermediate value theorem, there exists $q^1_1 \in [0.9345,0.935]$ such that $ \bar{h}^L_1(q^1_1)  = 0$. 
\item $x_1^1=0.939 ,x_1^2=0.94$ and $x_2^1= 0.959, x_2^2=0.961.$
For $i=1,2,$ it holds that
\begin{eqnarray*}
\bar{h}^U_2( x^i_1,x^1_2) \leq 0, \bar{h}^U_2(x^i_1,x^2_2) \geq 0, 
\end{eqnarray*}
and
\begin{eqnarray*}
 \bar{h}^L_2(x^1_1,x^i_2)  \leq 0, \bar{h}^L_2(x^2_1,x^i_2) \geq 0.
\end{eqnarray*}
By intermediate value theorem, there exists $(q^2_1,q^2_2) \in \big[ x^1_1,x^2_1 \big] \times \big[ x^1_2,x^2_2 \big]$ such that 
$\left \{ \begin{array}{lcl}
 \bar{h}^L_2\big (q^2_1,q^2_2\big) = 0,\\
 \bar{h}^U_2\big (q^2_1,q^2_2\big) = 0.
 \end{array} \right.$

\item Set $x_0=0.97$, then $\bar{h}^U_3(x_0) \leq 0$, where $(x_0,1)$ doesn't satisfy Condition $\varkappa(1).$ Also set $x'_0=0.972$, then $\bar{h}^U_3(x'_0) \geq 0$, where $(x'_0,1)$ satisfies Condition $\varkappa(1).$ Moreover, by intermediate value theorem, there exists $q^3_2 \in [0.97,0.972]$ such that $\bar{h}^U_3(q^3_2)= 0$. 
\end{enumerate}
Since $q^1_1 \leq 0.935<q^2_1<q^2_2<0.97\leq q^3_2,$ then the relation $ \mathcal{H}^1 \tilde{<}  \mathcal{H}^1 \tilde{<} \mathcal{H}^1$ holds. Then by Corollary \ref{nFRSB}, the Parisi measure is 3-FRSB with $\mathcal{F}= \{ 1,2  \}$(Figure \ref{3frsb12}). In this case, supp $\nu_P=\big \{0,1\big \} \cup \big[q^1_1,q^2_1 \big] \cup \big[q^2_2,q^3_2 \big].$

\end{example}

\begin{remark}
From the examples above, we notice that besides the 3-RSB phase, there are also three types of 3-FRSB appearing in the spherical model with 3 components. A smooth density can appear between the first and the second atoms of the Parisi measure. It can also appear between the second and the third atoms. The two smooth density can also appear simultaneously. 
It gives us a heuristic in mind that a smooth density may occur between any atoms of the Parisi measure at zero temperature. Also, multiple smooth density may occur among atoms simultaneously.   
We then expect that there will be $2^n-1$ types of $n$-FRSB appearing in the phases diagram of the spherical model with $n$ components.

\end{remark}

The rest of the paper is organized as follows: A general criterion for the Parisi measure $\nu_P$ to be $k-$RSB is derived in Section \ref{Seckrsb} while a general criterion for the Parisi measure $\nu_P$ to be $k-$FRSB is derived in Section \ref{Seckfrsb}.  A specific criterion for $\nu_P$ to be $n-$RSB is covered in Section \ref{Secnrsb}. Also, Corollary \ref{atmostnrsb} is proven at the end of Section \ref{Secnrsb}. Lastly, the proof of Theorem \ref{kRSB},\ref{kFRSB} and Corollary \ref{nRSB},\ref{nFRSB} are in Section \ref{Secproofmain}.

\section{Criteria for $k-$RSB Parisi measure}\label{Seckrsb}

Assume a measure $\nu$ is $k-$RSB and has the form of 
\begin{eqnarray*}
\nu(dx)=\sum^{k}_{i=1} m_i \cdot \mathbbm{1}_{[q_{i-1},q_i)}(x)dx+\Delta \cdot \delta_{\{1\}}(dx).
\end{eqnarray*}

Set $s=k$. For $\bar{q}=(0,q_1,q_2,\cdots,q_{k-1},1)$ and $l=1,\cdots,k-1$, we define
$$\dF^k_l=\dF^s_{l}(\bar{q}) \text{ and } \sA^k_l=\sA^s_{l}(\bar{q}).$$
We then define $$\mathcal{Z}^k=\mathcal{Z}^k(\bar{q}) \text{ and  }\mathcal{Y}^k=\mathcal{Y}^k(\bar{q}).$$
Moreover, set $s=k+1.$ For $l=1,\cdots,k$ and $\bar{x}_l=(0,q_1,\cdots,q_{l-1},x,q_l,\cdots,q_{k-1},1) \in [0,1]^{k+2}$ with $x \in (q_{l-1},q_l)$, we define
 \begin{eqnarray*}
 \sA_l^{k+1}(x)=\sA^{k+1}_{l}(0,q_1,\cdots,q_{l-1},x,q_l,\cdots,q_{k-1},1).
 \end{eqnarray*}
 and
  \begin{eqnarray*}
 \dF_l^{k+1}(x)=\dF^{k+1}_{l}(0,q_1,\cdots,q_{l-1},x,q_l,\cdots,q_{k-1},1).
 \end{eqnarray*}

We state the following criterion for the Parisi measure $\nu_P$ to be $k-$RSB:
\begin{theorem}\label{CriterionKRSB}
The Parisi measure $\nu_P$ for $\xi$ is $k-$RSB if and only if there exists unique $0<q_1<q_2<\cdots<q_{k-1}<1$ and $z_k>0$ such that the following conditions holds
\begin{enumerate}
\item for $l=1,\cdots,k,$ $
h\Big(q_{l-1},q_{l}, \big( \frac{(1+z_k)}{\dF^k_{l-\iota(k-l)}} \cdot r_2(q_{l-\iota(k-l-1)},q_{l-\iota(k-l)}) \big)^{(-1)^{k-l}}  \Big) = 0, 
$
\item $\mathcal{Y}^k \leq 1+z_k \leq \big(\mathcal{Z}^k\big)^{-1},$
\item for $l=1, \cdots, k$, $h\big(x^*_l,q_l,r_1(x^*_l,q_{l},q_{l-1}) \big)\geq0$
where  $x^*_l \in (q_{l-1},q_l)$  satisfies that $1+z_k=\big( \sA_l^{k+1}(x^*_l) \big)^{(-1)^{k-l}}.$

 \end{enumerate}

\end{theorem}

We first state a general characterization of Parisi measures:
\begin{theorem}[Theorem 2 in \cite{ArnabChen15}]\label{criterion}
For $\nu_P \in \mathcal{K},$ let 
\[
g(u)=\int^1_u \Big( \xi'(t)-\int^t_0 \frac{dr}{\nu_P((r,1])^2} \Big)dt.\]
 Then $\nu_P$ is the Parisi measure for the model $\xi$ if and only if the following conditions are satisfied:
 \begin{eqnarray*}
 &(i)& \quad \xi'(1)=\int^1_0 \frac{dr}{\nu_P((r,1])^2} , \\
 &(ii)& \text{ the function g satisfies } \min_{u \in [0,1]} g(u) \geq 0 \\
 &(iii)& \text{ for } S:= \{u \in [0,1):g(u)=0 \}, \text{ one has } \rho_P(S)=\rho_P([0,1)). \text{ Here } \rho_P \text{ is the measure} \\
 && \text{ induced by } \gamma_P, i.e.  \rho_P([0,s])=\gamma_P(s).
\end{eqnarray*}
\end{theorem}

We then come to the proof of Theorem \ref{CriterionKRSB} as follows:
\begin{proof}[Proof of Theorem \ref{CriterionKRSB}]

For $r \in [q_{l-1},q_l), l=1,2, \cdots ,k,$
\[ \nu((r,1])=m_l \cdot (q_l-r) + \sum^k_{i=l+1}m_i \cdot (q_i-q_{i-1})+ \Delta \]
and
\[ \int^{q_l}_{q_{l-1}}  \frac{dt}{\nu_P((t,1])^2} =\frac{(q_l-q_{l-1})}{ \big [\Delta+\sum^k_{j=l+1} m_j(q_j-q_{j-1}) \big] \big[ \Delta+\sum^k_{j=l} m_j(q_j-q_{j-1}) \big]} \]
Then a direct computation yields for $x \in [q_{l-1},q_{l}],$
\begin{eqnarray*}
g(x)&=& g(q_l) + \int^{q_l}_x  \bar{g}(t) dt \\
&=&\xi(q_l)-\xi(x)-(q_l-x) \cdot \sum^{l-1}_{i=1} \frac{(q_i-q_{i-1})}{\big [ \Delta+\sum^k_{j=i+1} m_j(q_j-q_{j-1}) \big] \big[\Delta +\sum^k_{j=i} m_j(q_j-q_{j-1}) \big]} \\
&&+\frac{q_l-x}{m_l \cdot \big[ \Delta+ \sum^k_{j=l}m_j(q_j-q_{j-1})   \big]} - \frac{1}{m^2_l} \ln \Big( 1+\frac{m_l (q_l-x)}{\Delta+ \sum^k_{j=l+1} m_j(q_j-q_{j-1}) } \Big) + g(q_l)
\end{eqnarray*}
where $\bar{g}(t)=\xi'(t)-\int^t_0 \frac{dr}{\nu_P((r,1])^2} $

By Theorem \ref{criterion}, $\nu$ is the Parisi measure for $\xi$ if and only if 
\begin{eqnarray}
&(1)& g(q_l)=\bar{g}(q_l)=0, \text{ for }l=1,\cdots k,  \label{Condition1} \\
&(2)& g(u) \geq 0, u\in [0,1].\label{Condition2}
\end{eqnarray}
Then by condition \eqref{Condition1},
\begin{eqnarray*}
 \bar{g}(q_l)-\bar{g}(q_{l-1})=0, \text{ for }l=1,\cdots k,
\end{eqnarray*}
which implies that 
\begin{eqnarray}\label{xip1}
\xi'(q_l)-\xi'(q_{l-1})= \frac{(q_l-q_{l-1})}{ \big [\Delta+\sum^k_{j=l+1} m_j(q_j-q_{j-1}) \big] \big[ \Delta+\sum^k_{j=l} m_j(q_j-q_{j-1}) \big]}, \text{ for }l=1,\cdots k, 
\end{eqnarray}

Set $z_l=\frac{m_l \cdot (q_l-q_{l-1})}{\Delta}$ for $l=1,2,\cdots,k.$
 \eqref{xip1} yields that 
\begin{eqnarray}\label{a1}
\Delta^2= \frac{(1-q_{k-1})}{[\xi'(1)-\xi'(q_{k-1})]} \cdot \frac{1}{(1+z_k)} 
\end{eqnarray}
and then for $l=1,\cdots,k-1,$
\begin{eqnarray}\label{a2}
1+\sum^k_{i=l}z_i&=&\frac{(q_l-q_{l-1})}{\Delta^2 \cdot [\xi'(q_l)-\xi'(q_{l-1})]} \cdot \frac{1}{(1+\sum^k_{i=l+1}z_i)} \nonumber \\
&=& \frac{[\xi'(1)-\xi'(q_{k-1})]\cdot (q_l-q_{l-1})}{ (1-q_{k-1})\cdot [\xi'(q_l)-\xi'(q_{l-1})]} \cdot \frac{(1+z_k)}{(1+\sum^k_{i=l+1}z_i)} ,
\end{eqnarray}

Based on the notations above and the equations \eqref{Condition1},\eqref{a1},\eqref{a2}, we rewrite $g(u)$ as follows, for $u \in [q_{l-1},q_l],l=1,2,\cdots, k,$
\begin{eqnarray*}
g(x)&=&\xi(q_l)-\xi(x)-\frac{(q_l-x)}{\Delta^2} \cdot \sum^{l-1}_{i=1} \frac{(q_i-q_{i-1})}{\big [ 1+\sum^k_{j=i+1} z_j \big] \big[1 +\sum^k_{j=i} z_j \big]} \\
&&+\frac{(q_l-q_{l-1})(q_l-x)}{\Delta^2 \cdot z_l  \big[ 1+ \sum^k_{j=l}z_j   \big]} - \frac{(q_l-q_{l-1})^2}{\Delta^2 \cdot z^2_l} \ln \Big( 1+\frac{(q_l-x)\cdot z_l}{ (q_l-q_{l-1})\big(1+ \sum^k_{j=l+1} z_j \big) } \Big) \\
&=&\xi(q_l)-\xi(x)-\xi'(q_{l-1})(q_l-x) + \frac{[\xi'(1)-\xi'(q_{k-1})](1+z_k)}{(1-q_{k-1})} \cdot \frac{(q_l-q_{l-1})(q_l-x)}{ z_l  \big[ 1+ \sum^k_{j=l}z_j   \big]}  \\
&&-  \frac{[\xi'(1)-\xi'(q_{k-1})]}{(1-q_{k-1})} \cdot \frac{(q_l-q_{l-1})^2(1+z_k)}{z^2_l} \ln \Big( 1+\frac{(q_l-x)\cdot z_l}{ (q_l-q_{l-1})\big(1+ \sum^k_{j=l+1} z_j \big) } \Big) 
\end{eqnarray*}
In particular, when $u=q_{l-1},$
\begin{eqnarray}
g(q_{l-1})&=&\xi(q_l)-\xi(q_{l-1})-\xi'(q_{l-1})(q_l-q_{l-1}) + \frac{[\xi'(1)-\xi'(q_{k-1})](1+z_k)}{(1-q_{k-1})} \cdot \frac{(q_l-q_{l-1})^2}{ z_l  \big[ 1+ \sum^k_{j=l}z_j   \big]} \nonumber \\
&&-  \frac{[\xi'(1)-\xi'(q_{k-1})]}{(1-q_{k-1})} \cdot \frac{(q_l-q_{l-1})^2(1+z_k)}{z^2_l} \ln \Big( 1+\frac{ z_l}{ 1+ \sum^k_{j=l+1} z_j  } \Big) =0.
\end{eqnarray}

Furthermore, by \eqref{a1},\eqref{a2}, it yields that
\begin{equation}\label{1+zl} 1+\sum^k_{i=l}z_i=  \left\{
\begin{array}{lcl}
&(1+z_k) \cdot  \prod^{\frac{k-l-2}{2}}_{i=0}\rr_{l+2i}(q),& \text{ for } k-l \text{  even, }  \\
&\prod^{\frac{k-l-1}{2}}_{i=0}\rr_{l+2i}(q),& \text{ for } k-l \text{  odd, } 
\end{array} \right. \end{equation}
and then
\begin{equation}\label{zl} z_l=  \left\{
\begin{array}{lcl}
&(1+z_k) \cdot  \prod^{\frac{k-l-2}{2}}_{i=0}\rr_{l+2i}(q)- \prod^{\frac{k-l-2}{2}}_{i=0}\rr_{l+1+2i}(q),& \text{ for } k-l \text{  even, }  \\
& \prod^{\frac{k-l-1}{2}}_{i=0}\rr_{l+2i}(q) -(1+z_k) \cdot  \prod^{\frac{k-l-3}{2}}_{i=0}\rr_{l+1+2i}(q),& \text{ for } k-l \text{  odd. } 
\end{array} \right. \end{equation}
 Notice that $\prod^{\frac{k-l-2}{2}}_{i=0}\rr_{l+2i}(q) \cdot \rr_{l+1+2i}(q)=\frac{(q_l-q_{l-1})[\xi'(q_k)-\xi'(q_{k-1})]}{[\xi'(q_l)-\xi'(q_{l-1})](q_k-q_{k-1})}$ for $k-l$ even and $\prod^{\frac{k-l-3}{2}}_{i=0}\rr_{l+1+2i}(q) \prod^{\frac{k-l-1}{2}}_{i=0}\rr_{l+2i}(q)=\frac{[\xi'(q_k)-\xi'(q_{k-1})](q_l-q_{l-1})}{(q_k-{q_{k-1}})[\xi'(q_l)-\xi'(q_{l-1})]}$ for $k-l$ odd.

Based on the notations above, we notice that
\begin{eqnarray}\label{1+zvs}
&&g(q_{l-1})\\
&=&
\left \{ \begin{array}{lcl}
h\Big(q_{l-1},q_{l}, (1+z_k)\prod^{\frac{k-l-2}{2}}_{i=0}\big(\frac{\rr_{l+2i}(q)}{\rr_{l+1+2i}(q)} \big) \Big) , \text{ if } k-l \text{ even}, \\
h \Big(q_{l-1},q_{l}, \frac{\rr_{k-1}(x)}{(1+z_k)} \prod^{\frac{k-l-3}{2}}_{i=0} \big( \frac{\rr_{l+2i}(q)}{\rr_{l+1+2i}(q)} \big)\Big) , \text{ if } k-l \text{ odd},
\end{array} \right. \nonumber \\
&=&
\left \{ \begin{array}{lcl}
h\Big(q_{l-1},q_{l}, \big(\frac{1+z_k}{\dF^k_l} \big) \cdot r_2(q_{l-1},q_l)  \Big) =h\Big(q_{l-1},q_{l}, \big( (1+z_k) \dF^k_{l-1} \cdot r_2(q_{l},q_{l-1})^{-1}  \Big) , \text{ if } k-l \text{ even}, \\ \
h\Big(q_{l-1},q_{l}, \big(\frac{\dF^k_{l-1}} {1+z_k}\big)  \cdot  r_2(q_{l},q_{l-1})^{-1} \Big)=h\Big(q_{l-1},q_{l}, \big(\dF^k_{l}(1+z_k)\big) ^{-1}  \cdot  r_2(q_{l-1},q_{l})\Big) , \text{ if } k-l \text{ odd}.
\end{array} \right. \nonumber
\end{eqnarray}

For $l = 1 \cdots k$, we define $h^k_l(\bar{q},1+z_k):=g(q_{l-1})$ to reflect the relation to $q_1,\cdots,q_k,1+z_k$ and $k,$ where $\bar{q}=(q_1,q_2,\cdots,q_{k-1}).$

Notice that $h(x,y,z)$ is strictly increasing with respect to $z$ when $z>0,$ we then have the following proposition:
\begin{lemma}\label{proph}
 $h^k_l(\bar{q},1+z_k)$ is strictly increasing with respect to $1+z_k$ if $k-l$ even. Similarly, $h^k_l(\bar{q},1+z_k)$ is strictly decreasing with respect to $1+z_k$ if $k-l$ is odd.
\end{lemma}

Since $\nu_P$ is a measure on $[0,1]$, it holds that $m_{l-1}<m_l<m_{l+1}$ for $l=1,\cdots,k-1$, which implies that
\begin{eqnarray*}
\frac{z_{l-1}}{(q_{l-1}-q_{l-2})}<\frac{z_l}{(q_l-q_{l-1})}<\frac{z_{l+1}}{(q_{l+1}-q_l)}.
\end{eqnarray*}

Recall the definition of  \eqref{dF},\eqref{sA}. Set $s=k,l=0.$ and apply it to $(q_0,q_1,\cdots,q_{s-1},q_s) \in [0,1]^{s+1}$.
Then by \eqref{zl}, $\text{ for } k-l \text{ even, }$
\begin{eqnarray}
1+z_k>\frac{[\xi'(q_l)-\xi'(q_{l-2})](q_{l-1}-q_{l-2})}{(q_l-q_{l-2})[\xi'(q_{l-1})-\xi'(q_{l-2})]} \cdot \prod^{\frac{k-l-2}{2}}_{i=0}\frac{\rr_{l+1+2i}(q)}{\rr_{l+2i}(q)}=\mathcal{A}^k_{l-1}(q),
\end{eqnarray}
and 
\begin{eqnarray}
1+z_k<\frac{[\xi'(q_{l+1})-\xi'(q_{l})](q_{l+1}-q_{l-1})}{(q_{l+1}-q_{l})[\xi'(q_{l+1})-\xi'(q_{l-1})]} \cdot \prod^{\frac{k-l-2}{2}}_{i=0}\frac{\rr_{l+1+2i}(q)}{\rr_{l+2i}(q)}=\frac{1}{\mathcal{A}^k_l(q)}.
\end{eqnarray}
where $q=(0,q_1,\cdots,q_{k-1},1).$

By \eqref{Condition2}, it holds that $g(u)\geq0$ for $u \in [q_{l-1},q_l]$ and $g(q_l)=0$, $l=1,2,\cdots, k$, which implies that $g''(q_l)\geq 0$ for $l=1,\cdots,k-1.$ We then obtain that
\begin{eqnarray*}
&&\xi''(q_l)\leq \frac{1}{ \big[\sum^k_{i=l+1}m_i \cdot (q_i-q_{i-1})+ \Delta \big]^2},\\
&<=>&\frac{[\xi'(1)-\xi'(q_{k-1})]}{\xi''(q_l) \cdot(1-q_{k-1})} \cdot (1+z_k)  \geq \Big[1+\sum^k_{i=l+1} z_i \Big]^2.
\end{eqnarray*}
By \eqref{1+zl}, for $k-l$ even,
\begin{eqnarray*}
1+z_k \geq \frac{\xi''(q_l)\cdot(q_k-q_{k-1})}{[\xi'(q_k)-\xi'(q_{k-1})]} \prod^{\frac{k-l-2}{2}}_{i=0} \rr^2_{l+1+2i}(q)=\dF^k_l(q),
\end{eqnarray*}
and
\begin{eqnarray*}
1+z_k \leq \frac{[\xi'(q_k)-\xi'(q_{k-1})]}{\xi''(q_{l-1})\cdot(q_k-q_{k-1})} \prod^{\frac{k-l-2}{2}}_{i=0} \frac{1}{\rr^2_{l+2i}(q)}=\frac{1}{\dF^k_{l-1}(q)}.
\end{eqnarray*}

Thus by Theorem \ref{criterion}, $\nu_P$ is the Parisi measure for $\nu_P$ if and only if  
\begin{enumerate}
\item $h^k_l(\bar{q}, 1+z_k)=0$ for $l=1,\cdots, k$
\item $\mathcal{Y}^k \leq 1+z_k \leq \big(\mathcal{Z}^k\big)^{-1},$
\item for $l=1, \cdots, k$, $g(x)>0$ for $x \in [q_{l-1},q_l].$

\end{enumerate}

It then remains for us to show that the condition (3) above is equivalent to the condition (3) in Theorem \ref{CriterionKRSB}. It is equivalent for us to show that $x^*_l \in (q_{l-1},q_l)$  satisfying $1+z_k=\sA^{k+1}_{l}(x^*_{l})$ are critical points of $g(x)$ and $g(x_l)=h^{k+1}_{l+1}(\bar{x}_l,\sA^{k+1}_l(x^*_l)),$ where $\bar{x}_l=(0,q_1,\cdots,q_{l-1},x^*_l,q_l,\cdots,q_{k-1},1).$

For $k-l$ even and the interval $[q_{l-1},q_l]$, a critical point $x^*_l \in [q_{l-1},q_l]$ of $g(x)$ satisfies that
\begin{eqnarray*}
&&-(1+z_k)^2[\xi'(x)-\xi'(q_{l-1})](q_l-x)\prod^{\frac{k-l-2}{2}}_{i=0}\rr_{l+2i}^2(q) -(x-q_{l-1})[\xi'(q_l)-\xi'(x)] \prod^{\frac{k-l-2}{2}}_{i=0}\rr_{l+1+2i}^2(q)\\
&&+(1+z_k)\frac{[\xi'(q_k)-\xi'(q_{k-1})](q_l-q_{l-1})}{(1-q_{k-1})[\xi'(q_l)-\xi'(q_{l-1})]}\Big[ (x-q_{l-1})[\xi'(q_l)-\xi'(x)]+[\xi'(x)-\xi'(q_{l-1})](q_l-x) \Big]=0
\end{eqnarray*}
and then
\begin{eqnarray}\label{1+zeven}
 1+z_k=\frac{(x^*_l-q_{l-1})[\xi'(q_l)-\xi'(x^*_l)]}{[\xi'(x^*_l)-\xi'(q_{l-1})](q_l-x^*_l)} \cdot \prod^{\frac{k-l-2}{2}}_{i=0}\frac{\rr_{l+1+2i}}{\rr_{l+2i}}. 
\end{eqnarray}
Notice that $1+z_k=\sA_l^{k+1}(x^*_l).$

For $k-l$ odd and the interval $[q_{l-1},q_{l}]$, a critical point $x^*_{l} \in [q_{l-1},q_{l}]$ of $g(x)$ satisfies that
\begin{eqnarray*}
&&-(1+z_k)^2[\xi'(q_{l})-\xi'(x)](x-q_{l-1})\prod^{\frac{k-l-1}{2}}_{i=1} \rr_{l-1+2i}^2 -[\xi'(x)-\xi'(q_{l-1})] (q_{l}-x)\prod^{\frac{k-l-1}{2}}_{i=0}\rr_{l+2i}^2\\
&&+(1+z_k)\frac{[\xi'(1)-\xi'(q_{k-1})](q_{l}-q_{l-1})}{(1-q_{k-1})[\xi'(q_{l})-\xi'(q_{l-1})]}\Big[ (x-q_{l-1})[\xi'(q_{l})-\xi'(x)]+[\xi'(x)-\xi'(q_{l-1})](q_{l}-x) \Big]=0
\end{eqnarray*}
and then
\begin{eqnarray}\label{1+zodd}
 1+z_k=\frac{[\xi'(x^*_{l})-\xi'(q_{l-1})](q_{l}-x^*_{l})}{(x^*_{l}-q_{l-1})[\xi'(q_{l})-\xi'(x^*_{l})]} \cdot \rr_{l-1} \cdot \prod^{\frac{k-l-1}{2}}_{i=0}\frac{\rr_{l+2i}}{\rr_{l-1+2i}}. 
\end{eqnarray}
Notice that $1+z_k=\big(\sA_l^{k+1}(x^*_l)\big)^{-1}.$

Now we substitue $1+z_k$ by \eqref{1+zeven} and \eqref{1+zodd} respectively in $g(x^*_l)$ and obtain that:
 \begin{eqnarray*}
&&g(x^*_l)\\
&=&\xi(q_l)-\xi(x^*_l)-\xi'(x^*_l)(q_l-x^*_l)+\frac{\big[\xi'(q_l)-\xi'(x^*_l)\big] \big[ \xi'(x^*_l)-\xi'(q_{l-1})\big] (q_l-x^*_l)(q_l-q_{l-1})}{\big \{ \big[\xi'(q_{l})-\xi'(q_{l-1}) \big](x^*_l-q_{l-1})-(q_{l}-q_{l-1}) \big[ \xi'(x^*_l)-\xi'(q_{l-1}) \big] \big \} } \\
&&-\frac{\big[ \xi'(q_{l})-\xi'(q_{l-1})\big] \big[\xi'(q_l)-\xi'(x^*_l)\big] \big[ \xi'(x^*_l)-\xi'(q_{l-1})\big] (q_l-x^*_l)(q_{l}-q_{l-1})(x^*_l-q_{l-1})}{\big \{  \big[\xi'(q_{l})-\xi'(q_{l-1}) \big](x^*_l-q_{l-1})-(q_{l}-q_{l-1}) \big[ \xi'(x^*_l)-\xi'(q_{l-1}) \big] \big \}^2 } \\
&&\cdot  \log \Big( \frac{\big[ \xi'(q_l)-\xi'(q_{l-1}) \big] (x^*_l-q_{l-1})}{ (q_l-q_{l-1})\big[\xi'(x^*_l)-\xi'(q_{l-1}) \big]} \Big)\\
&=&h\big(x^*_l,q_l,r_1(x^*_l,q_l ,q_{l-1})\big)\\
&=&h^{k+1}_{l+1}(\bar{x}_l,\sA^{k+1}_l(x_l)).
\end{eqnarray*}
which finishes the proof.

\end{proof}

\section{Criteria for $k-$FRSB Parisi measure}\label{Seckfrsb}

Assume a measure $\nu$ is $k-$FRSB and has the form of 
\begin{eqnarray*}
\nu(dx)=\sum^{t}_{j=1}\sum^{w_j}_{l=w_{j-1}+1} m_l \cdot \mathbbm{1}_{[q^j_{l-1},q^j_l)}(x)dx+\sum^{t-1}_{j=1} \omega(x) \cdot \mathbbm{1}_{[q^j_{w_j},q^{j+1}_{w_j})}(x)dx+\Delta \cdot \delta_{\{1\}}(dx).
\end{eqnarray*}

 For $j=1,\cdots,t,$ define $\bar{q}^j=(q^j_{w_{j-1}},q^j_{w_{j-1}+1},\cdots, q^j_{w_j})$,
we define, 
\begin{eqnarray*}
\nF^{j}_{l}=\dF^{s_j}_{l-w_{j-1}}(\bar{q}^j), \text{ for } w_{j-1} \leq l \leq w_j,
\end{eqnarray*}
and 
\begin{eqnarray*}
\nA^{j}_l =\sA^{s_j}_{l-w_{j-1}}(\bar{q}^j),\text{ for } w_{j-1}+1 \leq l \leq w_j-1.
\end{eqnarray*}
We then define
$$\nZ^{j}= \mathcal{Z}^{s_j}(q^j_{w_{j-1}},q^j_{w_{j-1}+1},\cdots, q^j_{w_j}),$$
and
$$\nY^{j}= \mathcal{Y}^{s_j}(q^j_{w_{j-1}},q^j_{w_{j-1}+1},\cdots, q^j_{w_j}).$$
Moreover, set $s=s_j+1.$ For $l=1,\cdots,s_j$ and $\bar{x}^j=(q^j_{w_{j-1}},q^j_{w_{j-1}+1},\cdots,q^j_{l-1},x,q^j_l,\cdots,q^j_{w_j}) \in [0,1]^{s_j+2}$ with $x \in (q^j_{l-1},q^j_l)$, we define
 \begin{eqnarray*}
 \nA_l^{j}(x)=\sA^{s_j+1}_{l}(\bar{x}^j).
 \end{eqnarray*}

We derive the following criterion for the Parisi measure $\nu_P$ to be $k-$FRSB:
\begin{theorem}\label{CriterionKFRSB}
The Parisi measure $\nu_P$ for $\xi$ is $k-$FRSB if and only if for $j=1,\cdots,t,$ $q^j_{w_{j-1}},q^j_{w_{j-1}+1},\cdots, q^j_{w_j}$ satisfy the following conditions
\begin{enumerate}
\item the system of equations
\begin{eqnarray} \label{conditionkfrsb}
\left \{ \begin{array}{lcl}
h\Big(q^j_{l-1},q^j_{l},  r_2(q^j_{l-1},q^j_l) (\nF^j_l)^{-1}  \cdot  {(*_j)^{(-1)^{w_j-l}}} \Big) = 0, \text{ for } l=1,\cdots,w_j, \\
\nF^{j}_{w_{j-1}}=\big(\nF^{j}_{w_j}\big)^{(-1)^{s_j}}, \text{ if } 1<j<t.
\end{array} \right.
\end{eqnarray}
where $*_j=\nF^j_{w_j}$ if $j=1,\cdots,t-1$ and $*_t=(\nF^t_{w_{t-1}})^{(-1)^{s_t}}.$
\item  $ \nY^j \leq *_j \leq (\nZ^j)^{-1} $.
\item for $l=w_{j_1}+1, \cdots, w_j$, $h\big(y^*_l,q^j_l,r_1(q^j_{l-1},q^j_l,y^*_l)\big)\geq0$
where  $y^*_l \in (q^j_{l-1},q^j_l)$  satisfies that $(*_j )^{(-1)^{\iota(w_j-l)}}=\nA_l^{j}(y^*_l).$
\item for $x \in [q^j_{w_j},q^{j+1}_{w_j}],$ $\frac{d^2}{dx^2} \big\{ \xi''(x)^{-\frac{1}{2}} \big \} \leq0. $

\end{enumerate}

\end{theorem}

\begin{proof}[Proof of Theorem \ref{CriterionKFRSB}]
Note that  for $ j_0=1,2, \cdots ,t+1,$ and $l_0=w_{j_0-1}+1, \cdots,w_{j_0},$ $r \in [q^{j_0}_{l_0-1},q^{j_0}_{l_0}),$
\[ \nu((r,1])=m_{l_0}  (q^{j_0}_{l_0}-r) + \sum^{w_{j_0}}_{l=l_0+1}m_l  (q^{j_0}_l-q^{j_0}_{l-1})+\sum^{t+1}_{j=j_0+1}\sum^{w_j}_{l=w_{j-1}+1}m_l(q^j_l-q^j_{l-1})+\sum^t_{j=j_0} \int^{q^{j+1}_{w_j}}_{q^j_{w_j}} \omega(s) ds + \Delta, \]
and
\begin{eqnarray}\label{re1}
 &&\int^{q^{j_0}_{l_0}}_{q^{j_0}_{l_0-1}}  \frac{dt}{\nu_P((t,1])^2} \nonumber \\
 &&=\frac{(q^{j_0}_{l_0}-q^{j_0}_{l_0-1})}{ \Big [ \sum^{w_{j_0}}_{l=l_0+1}m_l  (q^{j_0}_l-q^{j_0}_{l-1})+\sum^{t+1}_{j=j_0+1}\sum^{w_j}_{l=w_{j-1}+1}m_l(q^j_l-q^j_{l-1})+\sum^t_{j=j_0} \int^{q^{j+1}_{w_j}}_{q^j_{w_j}} \omega(s) ds + \Delta \Big]} \nonumber \\
&&\cdot \frac{1}{ \Big[ \sum^{w_{j_0}}_{l=l_0+1}m_l  (q^{j_0}_l-q^{j_0}_{l-1})+\sum^{t+1}_{j=j_0+1}\sum^{w_j}_{l=w_{j-1}+1}m_l(q^j_l-q^j_{l-1})+\sum^t_{j=j_0} \int^{q^{j+1}_{w_j}}_{q^j_{w_j}} \omega(s) ds + \Delta \Big]}.
\end{eqnarray}

By Theorem \ref{criterion}, $\nu_P$ is the Parisi measure for $\xi$ if and only if 
\begin{eqnarray}
&(1)& g(q^j_l)=\bar{g}(q^j_l)=0, \text{ for }j=1,\cdots t \text{ and } l=w_{j-1}+1, \cdots, w_j ,   \label{Condition1F}\\
&(2)& g(x)=\bar{g}(x)=\bar{g}'(x)=0, \text{ for } x \in [q^j_{w_j},q^{j+1}_{w_j}], j=1,\cdots t-1, \label{Condition2F}\\
&(3)& g(u) \geq 0, u\in [0,1].\label{Condition3F}\\
&(4)& \omega(x) \text{ is increasing in }[q^j_{w_j},q^{j+1}_{w_j}], \text{ for } j=1,\cdots,t-1.\label{Condition4F}
\end{eqnarray}
Then by condition \eqref{Condition2F}, $\text{ for } x \in [q^{j_0-1}_{w_{j_0-1}},q^{j_0}_{w_{j_0-1}}], j_0=1,\cdots t,$
\begin{eqnarray*}
\xi''(x)^{-\frac{1}{2}}=\nu((x,1])= \Delta+\sum^{t+1}_{j=j_0} \sum^{w_j}_{i=w_{j-1}+1}m_i (q^j_i-q^j_{i-1})+\sum^t_{j=j_0} \int^{q^{j+1}_{w_j}}_{q^j_{w_j}} \omega(s) ds +\int^{q^{j_0}_{w_{j_0-1}}}_x \omega(s)ds
\end{eqnarray*}
which implies that
\begin{eqnarray}\label{CF1}
\xi''(q^{j_0}_{w_{j_0-1}})^{-\frac{1}{2}}= \Delta+\sum^{t+1}_{j=j_0} \sum^{w_j}_{i=w_{j-1}+1}m_i (q^j_i-q^j_{i-1})+\sum^t_{j=j_0} \int^{q^{j+1}_{w_j}}_{q^j_{w_j}} \omega(s) ds 
\end{eqnarray}
and
\begin{eqnarray}\label{CF2}
\xi''(q^{j_0-1}_{w_{j_0-1}})^{-\frac{1}{2}}= \Delta+\sum^{t+1}_{j=j_0} \sum^{w_j}_{i=w_{j-1}+1}m_i (q^j_i-q^j_{i-1})+\sum^t_{j=j_0-1} \int^{q^{j+1}_{w_j}}_{q^j_{w_j}} \omega(s) ds.
\end{eqnarray}
Also we obtain that $\text{for } x \in [q^{j-1}_{w_{j-1}},q^{j}_{w_{j-1}}], j=1,\cdots t,$
$$\omega(x)=\frac{d}{dx} \big\{ \xi''(x)^{-\frac{1}{2}}\big \}=\frac{1}{2}\xi''(x)^{-\frac{3}{2}}\xi'''(x).$$
Thus condition \eqref{Condition4F} is equivalent to condition (4) in Theorem \ref{CriterionKFRSB}.

By subtracting \eqref{CF2} from \eqref{CF1}, we obtain that for $j=1,\cdots,t-1$, 
\begin{eqnarray*}
\xi''(q^{j}_{w_{j-1}})^{-\frac{1}{2}}-\xi''(q^{j}_{w_j})^{-\frac{1}{2}}=\sum^{w_j}_{i=w_{j-1}+1}m_i (q^j_i-q^j_{i-1}),
\end{eqnarray*}
and 
\begin{eqnarray*}
\xi''(q^{t}_{w_{t-1}})^{-\frac{1}{2}}= \Delta+ \sum^{w_{t}}_{i=w_{t-1}+1}m_i (q^{t}_i-q^{t}_{i-1}).
\end{eqnarray*}

Set $\rr^j_l=\frac{[\xi'(q^j_{l+1})-\xi'(q^j_l)](q^j_l-q^j_{l-1})}{(q^j_{l+1}-q^j_l)[\xi'(q^j_l)-\xi'(q^j_{l-1})]}$ . By induction, for $l=w_{t-1}+1,\cdots, w_{t},$
\begin{equation} \label{FR1}\Delta+\sum^{w_{t}}_{i=l} m_i(q^{t}_i-q^{t}_{i-1})= \left\{
\begin{array}{lcl}
&\frac{1}{\Big(\prod^{\frac{l-w_{t-1}-3}{2}}_{i=0}\rr^{t}_{w_{t-1}+1+2i}\Big)} \cdot \xi''(q^{t}_{w_{t-1}})^{-\frac{1}{2}},& \text{ for } l-w_{t-1} \text{  odd, }  \\
&\frac{(q^{t}_{l}-q^{t}_{l-1})}{\big[\xi'(q^{t}_{l})-\xi'(q^{t}_{l-1})\big]} \Big(\prod^{\frac{l-w_{t-1}-2}{2}}_{i=0}\rr^{t}_{w_{t-1}+1+2i} \Big)\cdot \xi''(q^{t}_{w_{t-1}})^{\frac{1}{2}},& \text{ for } l-w_{t-1} \text{  even. } 
\end{array} \right. \end{equation}


Also for $j=2,\cdots, t-1$ and $l=w_{j-1}+1,\cdots, w_j,w_j+1$, $\xi''(q^j_{w_j})^{-\frac{1}{2}}+\sum^{w_j}_{i=l} m_i(q^{j}_i-q^{j}_{i-1})$ can be expressed as
\begin{equation}\label{FR2}
\left\{
\begin{array}{lcl}
&\frac{1}{\Big(\prod^{\frac{l-w_{j-1}-3}{2}}_{i=0}\rr^{j}_{w_{j-1}+1+2i}\Big)} \cdot \xi''(q^{j}_{w_{j-1}})^{-\frac{1}{2}},& \text{ for } l-w_{j-1} \text{  odd, }  \\
&\frac{(q^{j}_{l}-q^{j}_{l-1})}{\big[\xi'(q^{j}_{l})-\xi'(q^{j}_{l-1})\big]} \Big(\prod^{\frac{l-w_{j-1}-2}{2}}_{i=0}\rr^{j}_{w_{j-1}+1+2i} \Big)\cdot \xi''(q^{j}_{w_{j-1}})^{\frac{1}{2}},& \text{ for } l-w_{j-1} \text{  even. } 
\end{array} \right. \end{equation}
In particular, for $j=2,\cdots,t-1$
\begin{equation} \label{FR3}\xi''(q^j_{w_j})^{-\frac{1}{2}}=
\left\{
\begin{array}{lcl}
&\frac{1}{\Big(\prod^{\frac{w_j-w_{j-1}-2}{2}}_{i=0}\rr^{j}_{w_{j-1}+1+2i}\Big)} \cdot \xi''(q^{j}_{w_{j-1}})^{-\frac{1}{2}},& \text{ for } s_j \text{  even, }  \\
&\frac{(q^{j}_{w_j}-q^{j}_{w_j-1})}{\big[\xi'(q^{j}_{w_j})-\xi'(q^{j}_{w_j-1})\big]} \Big(\prod^{\frac{w_j-w_{j-1}-3}{2}}_{i=0}\rr^{j}_{w_{j-1}+1+2i} \Big)\cdot \xi''(q^{j}_{w_{j-1}})^{\frac{1}{2}},& \text{ for } s_j \text{  odd. } 
\end{array} \right. \end{equation}
Notice that for $j=2,\cdots,t$, it then yields that $\nF^j_{w_{j-1}}=(\nF^j_{w_j})^{(-1)^{s_j}}.$

Based on the expressions above, we obtain that for $j=2,\cdots,t$ and $l=w_{j-1},\cdots,w_j$, $m_l$ is equal to
\begin{equation}\label{FR4}
\left\{
\begin{array}{lcl}
&\frac{\xi''(q^{j}_{w_{j-1}})^{-\frac{1}{2}}}{(q^j_l-q^j_{l-1})\Big(\prod^{\frac{l-w_{j-1}-3}{2}}_{i=0}\rr^{j}_{w_{j-1}+1+2i}\Big)} -\frac{\xi''(q^j_{w_{j-1}})^{\frac{1}{2}}}{\big[\xi'(q^{j}_{l})-\xi'(q^{j}_{l-1})\big]} \Big(\prod^{\frac{l-w_{j-1}-3}{2}}_{i=0}\rr^{j}_{w_{j-1}+1+2i} \Big) ,&  l-w_{j-1} \text{  odd, }  \\
&\frac{\xi''(q^{j}_{w_{j-1}})^{\frac{1}{2}}}{\big[\xi'(q^{j}_{l})-\xi'(q^{j}_{l-1})\big]} \Big(\prod^{\frac{l-w_{j-1}-2}{2}}_{i=0}\rr^{j}_{w_{j-1}+1+2i} \Big) -\frac{\xi''(q^{j}_{w_{j-1}})^{-\frac{1}{2}}}{(q^j_l-q^j_{l-1})\Big(\prod^{\frac{l-w_{j-1}-2}{2}}_{i=0}\rr^{j}_{w_{j-1}+1+2i}\Big)} ,& l-w_{j-1} \text{  even. } 
\end{array} \right. \end{equation}



Therefore, $\text{ for }j=1,\cdots t-1 \text{ and } l=w_{j-1}+1, \cdots, w_j $, \eqref{re1} can be simplified as
\begin{eqnarray*}
 \int^{q^{j}_{l}}_{q^{j}_{l-1}}  \frac{dt}{\nu_P((t,1])^2} =\frac{(q^{j}_{l}-q^{j}_{l-1})}{ \Big [ \sum^{w_j}_{i=l+1}m_i  (q^{j}_i-q^{j}_{i-1})+\xi''(q^{j}_{w_j})^{-\frac{1}{2}} \Big]\cdot \Big[ \sum^{w_j}_{i=l}m_i  (q^{j}_i-q^{j}_{i-1})+\xi''(q^{j}_{w_j})^{-\frac{1}{2}} \Big]}
\end{eqnarray*}
Also by Condition \eqref{Condition1F}, 
\begin{eqnarray*}
\bar{g}(q^j_{l})-\bar{g}(q^j_{l-1})=0,
\end{eqnarray*}
which implies that 
\begin{eqnarray}\label{xip}
\xi'(q^j_l)-\xi'(q^j_{l-1})= \frac{(q^{j}_{l}-q^{j}_{l-1})}{ \Big [ \sum^{w_j}_{i=l+1}m_i  (q^{j}_i-q^{j}_{i-1})+\xi''(q^{j}_{w_j})^{-\frac{1}{2}} \Big]\cdot \Big[ \sum^{w_j}_{i=l}m_i  (q^{j}_i-q^{j}_{i-1})+\xi''(q^{j}_{w_j})^{-\frac{1}{2}} \Big]}.
\end{eqnarray}

Similarly, for $l=w_{t-1}+1, \cdots, w_{t} $,
\begin{eqnarray*}
\xi'(q^{t}_l)-\xi'(q^{t}_{l-1})= \frac{(q^{t}_{l}-q^{t}_{l-1})}{ \Big [ \sum^{w_{t}}_{i=l+1}m_i  (q^{t}_i-q^{t}_{i-1})+\Delta \Big]\cdot \Big[ \sum^{w_{t}}_{i=l}m_i  (q^{t}_i-q^{t}_{i-1})+\Delta \Big]}.
\end{eqnarray*}

Apply \eqref{xip} for $j=1,\cdots,t-1$, and $l=w_{j-1},\cdots,w_j,$ we also obtain that
\begin{equation}\label{xi2}\xi''(q^j_{w_j})^{-\frac{1}{2}}+\sum^{w_j}_{i=l} m_i(q^{j}_i-q^{j}_{i-1})=
\left\{
\begin{array}{lcl}
&\Big( \prod^{\frac{w_j-l-1}{2}}_{i=0} \rr^j_{l+2i} \Big) \cdot \xi''(q^{j}_{w_{j}})^{-\frac{1}{2}},& \text{ for } w_j-l \text{  odd, }  \\
& \Big(\prod^{\frac{w_{j}-l-2}{2}}_{i=0}\rr^{j}_{l+2i} \Big) \cdot \frac{ (q^{j}_{w_j}-q^{j}_{w_j-1}) \cdot  \xi''(q^{j}_{w_{j}})^{\frac{1}{2}}  }{\big[\xi'(q^{j}_{w_j})-\xi'(q^{j}_{w_j-1})\big]} & \text{ for } w_{j}-l \text{  even. } 
\end{array} \right. \end{equation}
Then when $j=1$ and $l=1,\cdots,w_1,$ $m_l$ is also equal to
\begin{equation*}
\left\{
\begin{array}{lcl}
&\Big( \prod^{\frac{w_j-l-1}{2}}_{i=0} \rr^j_{l+2i} \Big) \cdot \frac{\xi''(q^{j}_{w_{j}})^{-\frac{1}{2}}}{(q^j_l-q^j_{l-1})} -  \Big(\prod^{\frac{w_{j}-l-1}{2}}_{i=0}\rr^{j}_{l+2i} \Big)^{-1} \cdot \frac{  \xi''(q^{j}_{w_{j}})^{\frac{1}{2}}  }{\big[\xi'(q^{j}_{l})-\xi'(q^{j}_{l-1})\big]} ,& \text{ for } w_j-l \text{  odd, }  \\
& \Big(\prod^{\frac{w_{j}-l-2}{2}}_{i=0}\rr^{j}_{l+1+2i} \Big)^{-1} \cdot \frac{   \xi''(q^{j}_{w_{j}})^{\frac{1}{2}}  }{\big[\xi'(q^{j}_{l})-\xi'(q^{j}_{l-1})\big]} -\Big( \prod^{\frac{w_j-l-2}{2}}_{i=0} \rr^j_{l+1+2i} \Big) \cdot \frac{\xi''(q^{j}_{w_{j}})^{-\frac{1}{2}}}{(q^j_l-q^j_{l-1})} & \text{ for } w_{j}-l \text{  even. } 
\end{array} \right. \end{equation*}

Then a direct computation yields for $j_0=1,\cdots t-1 , l=w_{j_0-1}+1, \cdots, w_{j_0} $ and $x \in [q^{j_0}_{{l_0}-1},q^{j_0}_{l_0}],$
\begin{eqnarray*}
&&g(x)\\
&&= g(q^{j_0}_{l_0}) + \int^{q^{j_0}_{l_0}}_x  \bar{g}(t) dt \\
&&=\xi(q^{j_0}_{l_0})-\xi(x) +\frac{q^{j_0}_{l_0}-x}{m_{l_0} \cdot \Big[ \sum^{w_j}_{i=l_0}m_i  (q^{j_0}_i-q^{j_0}_{i-1})+\xi''(q^{j_0}_{w_{j_0}})^{-\frac{1}{2}} \Big] } \\
&&-(q^{j_0}_{l_0}-x) \cdot \Big \{ \sum^{j_0-1}_{j=1} \sum^{w_j}_{l=w_{j-1}+1} \frac{(q^{j}_{l}-q^{j}_{l-1})}{ \Big [ \sum^{w_j}_{i=l+1}m_i  (q^{j}_i-q^{j}_{i-1})+\xi''(q^{j}_{w_j})^{-\frac{1}{2}} \Big]\cdot \Big[ \sum^{w_j}_{i=l}m_i  (q^{j}_i-q^{j}_{i-1})+\xi''(q^{j}_{w_j})^{-\frac{1}{2}} \Big]} \\
&&+\sum^{l_0-1}_{l=w_{j_0-1}+1} \frac{(q^{j_0}_l-q^{j_0}_{l-1})}{ \Big [ \sum^{w_{j_0}}_{i=l+1}m_i  (q^{j_0}_i-q^{j_0}_{i-1})+\xi''(q^{j_0}_{w_{j_0}})^{-\frac{1}{2}} \Big]\cdot \Big[ \sum^{w_{j_0}}_{i=l}m_i  (q^{j_0}_i-q^{j_0}_{i-1})+\xi''(q^{j_0}_{w_{j_0}})^{-\frac{1}{2}} \Big] } \Big \} \\
&&- \frac{1}{m^2_{l_0}} \ln \Big( 1+\frac{m_{l_0} (q^{j_0}_{l_0}-x)}{\xi''(q^{j_0}_{w_{j_0}})^{-\frac{1}{2}}+ \sum^{w_{j_0}}_{i=l_0+1} m_i(q^{j_0}_i-q^{j_0}_{j-1}) } \Big) 
\end{eqnarray*}
which implies that for $j=1,\cdots t-1 , l=w_{j-1}+1, \cdots, w_{j} $ and $x \in [q^{j}_{{l}-1},q^{j}_{l}],$
\begin{eqnarray*}
g(x)&=&\xi(q^{j}_{l})-\xi(x) -\xi'(q^{j}_{l-1})(q^{j}_{l}-x) +\frac{q^{j}_{l}-x}{m_{l} \cdot \Big[ \sum^{w_j}_{i=l}m_i  (q^{j}_i-q^{j}_{i-1})+\xi''(q^{j}_{w_j})^{-\frac{1}{2}} \Big] } \\
&&- \frac{1}{m^2_{l}} \ln \Big( 1+\frac{m_{l} (q^{j}_{l}-x)}{\xi''(q^{j}_{w_j})^{-\frac{1}{2}}+ \sum^{w_j}_{i=l+1} m_i(q^{j}_i-q^{j}_{i-1}) } \Big) .
\end{eqnarray*}
Here we use the relation $g(q^{j}_{l})=0.$

Similarly, for $l=w_{t-1}+1, \cdots, w_{t} $ and $x \in [q^{t}_{{l}-1},q^{t}_{l}],$
\begin{eqnarray*}
g(x)&=&\xi(q^{t}_{l})-\xi(x) -\xi'(q^{t}_{l-1})(q^{t}_{l}-x) +\frac{q^{t}_{l}-x}{m_{l} \cdot \Big[ \sum^{w_{t}}_{i=l}m_i  (q^{t}_i-q^{t}_{i-1})+\Delta \Big] } \\
&&- \frac{1}{m^2_{l}} \ln \Big( 1+\frac{m_{l} (q^{t}_{l}-x)}{\Delta+ \sum^{w_{t}}_{i=l+1} m_i(q^{t}_i-q^{t}_{i-1}) } \Big) .
\end{eqnarray*}

In particular, for $j=1,\cdots t-1 , l=w_{j-1}+1, \cdots, w_{j} $ and $x = q^{j}_{{l}-1},$
\begin{eqnarray*}
g(q^{j}_{{l}-1})&=&\xi(q^{j}_{l})-\xi(q^{j}_{{l}-1}) -\xi'(q^{j}_{l-1})(q^{j}_{l}-q^{j}_{{l}-1}) +\frac{q^{j}_{l}-q^{j}_{{l}-1}}{m_{l} \cdot \Big[ \sum^{w_j}_{i=l}m_i  (q^{j}_i-q^{j}_{i-1})+\xi''(q^{j}_{w_j})^{-\frac{1}{2}} \Big] } \\
&&- \frac{1}{m^2_{l}} \ln \Big( 1+\frac{m_{l} (q^{j}_{l}-q^{j}_{{l}-1})}{\xi''(q^{j}_{w_j})^{-\frac{1}{2}}+ \sum^{w_j}_{i=l+1} m_i(q^{j}_i-q^{j}_{i-1}) } \Big)=0 .
\end{eqnarray*}
Also for $l=w_{t-1}+1, \cdots, w_{t} $ and $x = q^{t}_{{l}-1},$
\begin{eqnarray*}
g( q^{t}_{{l}-1})&=&\xi(q^{t}_{l})-\xi( q^{t}_{{l}-1}) -\xi'(q^{t}_{l-1})(q^{t}_{l}- q^{t}_{{l}-1}) +\frac{q^{t}_{l}- q^{t}_{{l}-1}}{m_{l} \cdot \Big[ \sum^{w_{t}}_{i=l}m_i  (q^{t}_i-q^{t}_{i-1})+\Delta \Big] } \\
&&- \frac{1}{m^2_{l}} \ln \Big( 1+\frac{m_{l} (q^{t}_{l}- q^{t}_{{l}-1})}{\Delta+ \sum^{w_{t}}_{i=l+1} m_i(q^{t}_i-q^{t}_{i-1}) } \Big) =0.
\end{eqnarray*}

Based on the relations \eqref{FR1},\eqref{FR2},\eqref{FR3},\eqref{FR4}, we notice that for $j=2,\cdots,t$ and $l=w_{j-1}+1,\cdots,w_j$, it holds that $$g(q^j_{l-1})=h \Big(q^j_{l-1},q^j_{l}, \big(\nF^j_{{l}}\big)^{-1}r_{2}(q^j_{l-1},q^j_{l})  \cdot \big(\nF^j_{w_{j-1}} \big)^{(-1)^{s_j+\iota(w_j-l)}} \Big) =0,$$
and for $j=1$, $l=1,\cdots,w_1,$
\begin{eqnarray*}
g(q^j_{l-1})=h \Big(q^1_{l-1},q^1_l, (F^1_l)^{-1} r_2(q^1_{l-1},q^1_l ) \cdot  \big( \nF^1_{w_1} \big)^{(-1)^{w_1-l}} \Big)=0.
\end{eqnarray*}

For $s=s_j$ and $l = w_{j-1}+1, \cdots w_j$, we define $h^j_l(\bar{q}^j, F_{w_{j-1}}):=g(q^j_{l-1})$ to reflect the relation to $\bar{q}^j=(q^j_{w_{j-1}},q^j_{w_{j-1}+1},\cdots,q^j_{w_j})$ and $F_{w_{j-1}}.$

By \eqref{Condition3F}, it holds that $g(u)\geq0$ for $u \in [0,1]$ and $g(q^j_{l})=0$, $\text{ for }j=1,\cdots t , l=w_{j-1}, \cdots, w_j $, which implies that $g''(q^j_l)\geq 0$. We obtain that
\begin{eqnarray*}
\xi''(q^j_l)\leq \frac{1}{ \big[ \sum^{w_j}_{i=l+1}m_i  (q^{j}_i-q^{j}_{i-1})+\xi''(q^{j}_{w_j})^{-\frac{1}{2}} \big]^2}.
\end{eqnarray*}
By \eqref{xi2}, for $j=1,\cdots,t-1,$ it then yields that if $w_j-l$ odd,
\begin{eqnarray*}
 \nF^j_{w_{j}} \cdot \nF^j_{l}  \leq 1,
\end{eqnarray*}
and if $w_j-l$ even,
\begin{eqnarray*}
\nF^j_{l}  \leq  \nF^j_{w_{j}} .
\end{eqnarray*}

Since $\nu_P$ is a measure on $[0,1]$, it holds that $m_{l-1}<m_l<m_{l+1}$ for $l=1,\cdots,k-1$, which implies that
by \eqref{FR4}, for $j=1,\cdots,t-1,l=w_{j-1}+1,\cdots,w_j$, if $w_{j}-l$ even,
\begin{eqnarray}
 \nF^j_{w_{j}} \cdot \nA^j_{l} \leq 1,
\end{eqnarray}
and if $w_{j}-l$ odd,
\begin{eqnarray}
\nF^j_{w_{j}} \geq \nA^j_l.
\end{eqnarray}

Similarly for $j=t$ and $l=w_{t-1}+1,\cdots,w_t,$ it yields that if $k-l$ odd,
\begin{eqnarray*}
\nA^t_l \leq (\nF^t_{w_{t-1}})^{(-1)^{s_t}}  \leq (\nF^t_{l} )^{-1}, 
\end{eqnarray*}
and if $k-l$ even,
\begin{eqnarray}
\nF^t_{l} \leq (\nF^t_{w_{t-1}})^{(-1)^{s_t}} \leq (\nA^t_{l})^{-1} .
\end{eqnarray}

The inequalities above are equivalent to the relation $$\nY^j \leq *_j \leq (\nZ^j )^{-1}, \text{ for } j=1,\cdots,t.$$

Now we show that $g(x) \geq 0$ for $x \in (q^j_{l-1},q^j_l)$ is equivalent to $h\big(y^*_l,q^j_l,r_1(q^j_{l-1},q^j_l,y^*_l)\big)\geq0$, with  $y^*_l \in (q^j_{l-1},q^j_l)$  satisfying $F^j_{w_j}=\nA_l^{k+1}(y^*_l)$ for $j=1,\cdots,t$ and $l=w_{j-1}+1,\cdots,w_j.$

Fix $j=1,\cdots t-1 , l=w_{j-1}+1, \cdots, w_j $ and the interval $[q^{j}_{{l}-1},q^{j}_{l}].$
For $w_j-l$ even, a critical point $y^*_l \in [q^j_{l-1},q^j_l]$ of $g(x)$ satisfies that
\begin{eqnarray*}
&&-(\nF^j_{w_{j}})^{2}[\xi'(x)-\xi'(q^j_{l-1})](q^j_l-x)\prod^{\frac{w_j-l-2}{2}}_{i=0}(\rr_{l+2i}^j)^2 -(x-q^j_{l-1})[\xi'(q^j_l)-\xi'(x)] \prod^{\frac{w_j-l-2}{2}}_{i=0}(\rr_{l+1+2i}^j)^2 \\
&&+\nF^j_{w_{j}} \cdot \prod^{\frac{w_j-l-2}{2}}_{i=0} (\rr_{l+2i}^j \cdot \rr_{l+1+2i}^j) \Big[ (x-q^j_{l-1})[\xi'(q^j_l)-\xi'(x)]+[\xi'(x)-\xi'(q^j_{l-1})](q^j_l-x) \Big]=0,
\end{eqnarray*}
and then
\begin{eqnarray*}
F^j_{w_j}=\frac{\big[\xi'(q^j_l)-\xi'(y^*_l)\big](y^*_l-q^j_{l-1})}{(q^j_l-y^*_l) \big[\xi'(y^*_l)-\xi'(q^j_{l-1})]} \cdot \prod^{\frac{w_j-l-2}{2}}_{i=0} \frac{\rr^j_{l+1+2i}}{\rr^j_{l+2i}}=\nA_l^{j}(y^*_l). 
\end{eqnarray*}


For $w_j-l$ odd, a critical point $y^*_l \in [q^j_{l-1},q^j_l]$ of $g(x)$ satisfies that
\begin{eqnarray*}
&&-(\nF^j_{w_{j}})^{2}[\xi'(q^j_{l})-\xi'(x)](x-q^j_{l-1}) \prod^{\frac{w_j-l-3}{2}}_{i=0}(\rr_{l+1+2i}^j)^2  -[\xi'(x)-\xi'(q^j_{l-1})] (q^j_{l}-x) \prod^{\frac{w_j-l-1}{2}}_{i=0}(\rr_{l+2i}^j)^{2}\\
&& +\nF^j_{w_{j}} \cdot \prod^{\frac{w_j-l-3}{2}}_{i=0} \rr_{l+1+2i}^j \prod^{\frac{w_j-l-1}{2}}_{i=0} \rr_{l+2i}^j \cdot \Big[ (x-q^j_{l-1})[\xi'(q^j_{l})-\xi'(x)]+[\xi'(x)-\xi'(q^j_{l-1})](q^j_{l}-x) \Big]=0,
\end{eqnarray*}
and then
\begin{eqnarray*}
F^j_{w_j}=\frac{(q^j_l-y^*_l) \big[\xi'(y^*_l)-\xi'(q^j_{l-1})]}{\big[\xi'(q^j_l)-\xi'(y^*_l)\big](y^*_l-q^j_{l-1})} \cdot \frac{\prod^{\frac{w_j-l-1}{2}}_{i=0}  \rr^j_{l+2i}}{\prod^{\frac{w_j-l-3}{2}}_{i=0}  \rr^j_{l+1+2i}}=\big(\nA_l^{j}(y^*_l) \big)^{-1}. 
\end{eqnarray*}


Similarly, for $j=t , l=w_{t-1}+1, \cdots, w_t $ and the interval $[q^{t}_{{l}-1},q^{t}_{l}].$
For $k-l$ even, a critical point $y^*_l \in [q^t_{l-1},q^t_l]$ of $g(x)$ satisfies that
\begin{eqnarray*}
&&- \big(\nF^t_{w_{t-1}} \big)^{2(-1)^{s_t}}  [\xi'(x)-\xi'(q^t_{l-1})](q^t_l-x)\prod^{\frac{k-l-2}{2}}_{i=0}(\rr_{l+2i}^t)^2 -(x-q^t_{l-1})[\xi'(q^t_l)-\xi'(x)] \prod^{\frac{k-l-2}{2}}_{i=0}(\rr_{l+1+2i}^t)^2 \\
&&+\big(\nF^t_{w_{t-1}} \big)^{(-1)^{s_t}} \cdot \prod^{\frac{k-l-2}{2}}_{i=0} (\rr_{l+2i}^t \cdot \rr_{l+1+2i}^t) \Big[ (x-q^t_{l-1})[\xi'(q^t_l)-\xi'(x)]+[\xi'(x)-\xi'(q^t_{l-1})](q^t_l-x) \Big]=0,
\end{eqnarray*}
and then
\begin{eqnarray*}
\big(\nF^t_{w_{t-1}} \big)^{(-1)^{s_t}}=\frac{\big[\xi'(q^t_l)-\xi'(y^*_l)\big](y^*_l-q^t_{l-1})}{(q^t_l-y^*_l) \big[\xi'(y^*_l)-\xi'(q^t_{l-1})]} \cdot \prod^{\frac{k-l-2}{2}}_{i=0} \frac{\rr^t_{l+1+2i}}{\rr^j_{l+2i}}=\nA_l^{t}(y^*_l). 
\end{eqnarray*}
For $k-l$ odd, a critical point $y^*_l \in [q^t_{l-1},q^t_l]$ of $g(x)$ satisfies that
\begin{eqnarray*}
&&- \big(\nF^t_{w_{t-1}} \big)^{2(-1)^{s_t}} [\xi'(q^t_{l})-\xi'(x)](x-q^t_{l-1}) \prod^{\frac{k-l-3}{2}}_{i=0}(\rr_{l+1+2i}^t)^2  -[\xi'(x)-\xi'(q^t_{l-1})] (q^t_{l}-x) \prod^{\frac{k-l-1}{2}}_{i=0}(\rr_{l+2i}^t)^{2}\\
&& + \big(\nF^t_{w_{t-1}} \big)^{(-1)^{s_t}}   \prod^{\frac{k-l-3}{2}}_{i=0} \rr_{l+1+2i}^t \prod^{\frac{k-l-1}{2}}_{i=0} \rr_{l+2i}^t  \Big[ (x-q^t_{l-1})[\xi'(q^t_{l})-\xi'(x)]+[\xi'(x)-\xi'(q^t_{l-1})](q^t_{l}-x) \Big]=0,
\end{eqnarray*}
and then
\begin{eqnarray*}
 \big(\nF^t_{w_{t-1}} \big)^{(-1)^{s_t}} =\frac{(q^t_l-y^*_l) \big[\xi'(y^*_l)-\xi'(q^t_{l-1})]}{\big[\xi'(q^t_l)-\xi'(y^*_l)\big](y^*_l-q^t_{l-1})} \cdot \frac{\prod^{\frac{k-l-1}{2}}_{i=0}  \rr^t_{l+2i}}{\prod^{\frac{k-l-3}{2}}_{i=0}  \rr^t_{l+1+2i}}=\big(\nA_l^{t}(y^*_l) \big)^{-1}. 
\end{eqnarray*}

Now for $j=1,\cdots,t,$ we substitue $*_j$ by $\nA^{j}_l(y^*_l)$ in $g(y^*_l)$ and obtain that:
 \begin{eqnarray*}
&&g(y^*_l)\\
&=&\xi(q^j_l)-\xi(y^*_l)-\xi'(y^*_l)(q^j_l-y^*_l)+\frac{\big[\xi'(q^j_l)-\xi'(y^*_l)\big] \big[ \xi'(y^*_l)-\xi'(q^j_{l-1})\big] (q^j_l-y^*_l)(q^j_l-q^j_{l-1})}{\big \{ \big[\xi'(q_{l})-\xi'(q^j_{l-1}) \big](y^*_l-q^j_{l-1})-(q_{l}-q^j_{l-1}) \big[ \xi'(y^*_l)-\xi'(q^j_{l-1}) \big] \big \} } \\
&&-\frac{\big[ \xi'(q_{l})-\xi'(q^j_{l-1})\big] \big[\xi'(q^j_l)-\xi'(y^*_l)\big] \big[ \xi'(y^*_l)-\xi'(q^j_{l-1})\big] (q^j_l-y^*_l)(q_{l}-q^j_{l-1})(y^*_l-q^j_{l-1})}{\big \{  \big[\xi'(q_{l})-\xi'(q^j_{l-1}) \big](y^*_l-q^j_{l-1})-(q_{l}-q^j_{l-1}) \big[ \xi'(y^*_l)-\xi'(q^j_{l-1}) \big] \big \}^2 } \\
&&\cdot  \log \Big( \frac{\big[ \xi'(q^j_l)-\xi'(q^j_{l-1}) \big] (y^*_l-q^j_{l-1})}{ (q^j_l-q^j_{l-1})\big[\xi'(y^*_l)-\xi'(q^j_{l-1}) \big]} \Big)\\
&=&h\big(y^*_l,q^j_l,r_1(q^j_{l-1},q^j_l,y^*_l)\big).
\end{eqnarray*}
Therefore the condition $h\big(y^*_l,q^j_l,r_1(q^j_{l-1},q^j_l,y^*_l)\big)>0$ is equivalent to the condition $g(x)>0$ for $x \in [q^j_{l-1},q^j_l].$

\end{proof}

\section{A special case: $k=n$} \label{Secnrsb}
 
 Recall the setting of the model
 \begin{eqnarray}
\xi(x):=\sum^{n-1}_{l=1} \lambda_l \cdot x^{p_l} +\Big( 1- \sum^{n-1}_{l=1} \lambda_l \Big) \cdot x^{p_n}.
\end{eqnarray}
Then when $k=n,$ Theorem \ref{CriterionKRSB} can be simplified as follows:
\begin{theorem}\label{nrsb}
The Parisi measure $\nu_P$ for $\xi$ is $n-$RSB if and only if there exists unique increasing sequence $0<q_1<q_2<\cdots<q_{n-1}<1$ and $z_n>0$ such that the following conditions holds
\begin{enumerate}
\item  the system of equations
\begin{eqnarray*}\left\{
\begin{array}{lcl}
h^n_1(\bar{q}, 1+z_n)=0,\\
\qquad \cdots \cdots\\
h^n_n(\bar{q}, 1+z_n)=0,
\end{array} \right . 
\end{eqnarray*}
where $\bar{q}=(q_1,q_2,\cdots,q_{k-1}).$
\item $\mathcal{Y}^n \leq 1+z_n \leq \big(\mathcal{Z}^n\big)^{-1}$. \label{bound}
\end{enumerate}
\end{theorem}

\begin{proof}
By Theorem \ref{CriterionKRSB}, it suffices for us to show that condition (3) can be deduced from conditions (1) and (2).

First note that
\begin{eqnarray*}
\frac{d^2}{dx^2} \big \{ \xi''(x)^{-\frac{1}{2}} \big \}&=&\frac{1}{4} \xi''(x)^{-\frac{5}{2}} \cdot \big \{ 3 \xi'''(x)^2-2\xi''(x)\xi''''(x) \big\}\\
&=&\frac{1}{4} \xi''(x)^{-\frac{5}{2}} \cdot \Big (   \sum^n_{l=1}p_l^3(p_l-1)^2(p_l-2)\lambda^2_l x^{2p_l-6} \\
&&-2 \sum_{1 \leq i<j \leq n}p_i(p_i-1)p_j(p_j-1)(p^2_i+p^2_j+p_i+p_j-3p_ip_j) \lambda_i \lambda_j x^{p_i+p_j-6}\Big)
\end{eqnarray*}
Then by Descartes' rule of signs,  $\frac{d^2}{dx^2} \big \{ \xi''(x)^{-\frac{1}{2}} \big \}$ has at most $2(n-1)$ strictly positive roots.

Recall that 
 \begin{eqnarray*}
h(x,y,z)=\xi(y)-\xi(x)-\xi'(x)(y-x)+\big[\xi'(y)-\xi'(x)\big]  (y-x) \Big( \frac{1}{( z-1 ) } -\frac{z}{(z-1)^2 }  \log ( z ) \Big),
\end{eqnarray*}
and $ r_2(x,y)=\frac{\xi''(y)(y-x)}{[\xi'(y)-\xi'(x)]}.$ Define $h^F(x,y)=h(x,y,r_2(x,y)).$

By computation, it yields that 
\begin{eqnarray}\label{h^F(x,y)}
&&\text{sgn} \Big( \frac{\partial}{\partial x} h^F(x,y)  \Big) = \text{sgn} \Big( [\xi'(y)-\xi'(x)]^2-\xi''(x)\xi''(y)(y-x)^2 \Big) \nonumber \\
&&\text{sgn} \Big( \frac{\partial}{\partial y} h^F(x,y)  \Big) = \text{sgn} \Big( \xi'''(y)[\xi'(y)-\xi'(x)](y-x)+2\xi''(y)\big[ \xi'(y)-\xi'(x)-\xi''(y)(y-x) \big] \Big)  \nonumber
\end{eqnarray}
and
\begin{eqnarray*}
&&\frac{\partial^2}{\partial x^2} h^F(x,y) \big|_{x=y}=\frac{\partial^3}{\partial x^3} h^F(x,y) \big|_{x=y}=\frac{\partial^2}{\partial y^2} h^F(x,y) \big|_{y=x}=0,\\
&&\text{sgn} \Big( \frac{\partial^4}{\partial x^4} h^F(x,y) \big|_{x=y} \Big) = \text{sgn} \Big( 3 \xi'''(y)^2-2\xi''(y)\xi''''(y) \Big),\\
&&\text{sgn} \Big( \frac{\partial^3}{\partial y^3} h^F(x,y) \big|_{y=x} \Big) = -\text{sgn} \Big( 3 \xi'''(x)^2-2\xi''(x)\xi''''(x) \Big).
\end{eqnarray*}
If for fixed $x_0$, $y_0>x_0$ is a critical point of $h^F(x,y)$ with respect to $y$, then
\begin{eqnarray*}
&&\text{sgn} \Big( \frac{\partial^2}{\partial y^2}h^F(x_0,y) \big|_{y=y_0}\Big) \\
&=&\text{sgn} \Big( \xi''''(y_0) \big[ \xi'(y_0)-\xi'(x_0)\big] (y_0-x_0)+3\xi'''(y_0) \big[ \xi'(y_0)-\xi'(x_0)-\xi''(y_0)(y_0-x_0) \Big) \\
&=&\text{sgn} \Big( 2\xi''(y_0) \xi''''(y_0) -3\xi'''(y_0)^2 \Big).
\end{eqnarray*}
Here we use the relation $\xi'''(y_0)[\xi'(y_0)-\xi'(x_0)](y_0-x_0)=-2\xi''(y_0)\big[ \xi'(y_0)-\xi'(x_0)-\xi''(y_0)(y_0-x_0) \big]$.

Thus for fixed $x$, a local minimum $y>x$ of $h^F(x,y)$ satisfies $\frac{d^2}{dy^2} \big \{ \xi''(y)^{-\frac{1}{2}} \big \}<0$ while a local maximum $y>x$ of $h^F(x,y)$ satisfies $\frac{d^2}{dy^2} \big \{ \xi''(y)^{-\frac{1}{2}} \big \}>0.$

Also $h^n_{l}(\bar{q}, \dF^n_{l})$ and $h^n_{l}(\bar{q}, \dF^n_{l-1})$ can be computed as follows: 
\begin{eqnarray*}
h^n_{l}(\bar{q}, \dF^n_{l})&=&\xi(q_l)-\xi(q_{l-1})-\xi'(q_{l-1})(q_l-q_{l-1})+\frac{(q_l-q_{l-1})\big[\xi'(q_l)-\xi'(q_{l-1}) \big]^2}{\big \{\xi''(q_l)(q_l-q_{l-1})-[\xi'(q_l)-\xi'(q_{l-1})] \big \} }\\
&&-\frac{\xi''(q_l)[\xi'(q_l)-\xi'(q_{l-1})]^2(q_l-q_{l-1})^2}{\big \{\xi''(q_l)(q_l-q_{l-1})-[\xi'(q_l)-\xi'(q_{l-1})] \big \}^2 } \log\Big(\frac{\xi''(q_l)(q_l-q_{l-1})}{[\xi'(q_l)-\xi'(q_{l-1})]} \Big),
\end{eqnarray*}
and
\begin{eqnarray*}
h^n_{l}(\bar{q}, \mathcal{F}^n_{l-1})&=&\xi(q_l)-\xi(q_{l-1})-\xi'(q_{l-1})(q_l-q_{l-1})+\frac{\xi''(q_{l-1})\big[\xi'(q_l)-\xi'(q_{l-1}) \big](q_l-q_{l-1})^2}{\big \{\xi'(q_l)-\xi'(q_{l-1})-\xi''(q_{l-1})(q_l-q_{l-1}) \big \} }\\
&&-\frac{\xi''(q_{l-1})[\xi'(q_l)-\xi'(q_{l-1})]^2(q_l-q_{l-1})^2}{\big \{\xi'(q_l)-\xi'(q_{l-1})-\xi''(q_{l-1})(q_l-q_{l-1}) \big \}^2 } \log\Big(\frac{[\xi'(q_l)-\xi'(q_{l-1})]}{\xi''(q_{l-1})(q_l-q_{l-1})} \Big).
\end{eqnarray*}

Notice that $h^n_{l}(\bar{q}, \dF^n_{l})=h^F(q_{l-1},q_l)$ and $h^n_{l}(\bar{q}, \dF^n_{l-1})=-h^F(q_{l},q_{l-1})$. Also, 
\begin{eqnarray*}
 \text{sgn} \Big( \frac{\partial}{\partial q_{l-1}} h^n_l(\bar{q},F^n_l)  \Big)&=&-\text{sgn} \Big(  \frac{\partial}{\partial q_{l}}  h^n_{l}(\bar{q}, \dF^n_{l-1}) \Big),\\
 &=& 1-\dF^n_l\cdot \dF^n_{l-1}.
\end{eqnarray*}

Indeed, there must exist $x \in (q_{l-1},q_l)$ such that $\frac{d^2}{dx^2} \big \{ \xi''(x)^{-\frac{1}{2}} \big \}>0.$ Otherwise, if for all $x \in [q_{l-1},q_l),\frac{d^2}{dx^2} \big \{ \xi''(x)^{-\frac{1}{2}} \big \}<0,$ it yields that $-h^F(q_{l},x) <0$ and $h^F(x,q_l)>0$. Then by Lemma \ref{proph}, it holds that $\dF^n_{l-1}\cdot \dF^n_l>1$, which leads to a contradiction with $\mathcal{Y}^n  \leq \big(\mathcal{Z}^n\big)^{-1}$.

Now for any $0<x_1<\cdots<x_{n-1}<1$, set $\bar{x}=(0,x_1,x_2,\cdots,x_{n-1},1).$

For $l=1, h^n_{1}(\bar{x}, \dF^n_{1})$ is only with respect to $x_1.$ We then denote the first local minimum of $h^n_{1}(\bar{x}, \dF^n_{1})$ by $m_1$ and the first local maximum by $M_1.$  By \eqref{h^F(x,y)}, it holds that $m_1<M_1.$  We then obtain that 
\begin{eqnarray*}
\frac{d^2}{dx^2} \big \{ \xi''(x)^{-\frac{1}{2}} \big \} \big|_{x=m_1}<0 \text{ and } \frac{d^2}{dx^2} \big \{ \xi''(x)^{-\frac{1}{2}} \big \} \big|_{x=M_1}>0.
\end{eqnarray*}
Since $\lim_{x \rightarrow 0^+} \frac{d^2}{dx^2} \big \{ \xi''(x)^{-\frac{1}{2}} \big \} \big|_{x=m^1_1}=+\infty$, by intermediate value theorem, there exists zeros in intervals $[0,m_1]$ and $[m_1,M_1]$ respectively.

Now for $l=2,\cdots, n-1,$ denote the first local minimum larger than $M_{l-1}$ of $h_l(\bar{x},F_l)$ with respect to $x_l$ by $m_l$ and the first local maximum by $M_l.$ By \eqref{h^F(x,y)}, it holds that $m_l<M_l.$ We then obtain that 
\begin{eqnarray*}
\frac{d^2}{dx^2} \big \{ \xi''(x)^{-\frac{1}{2}} \big \} \big|_{x=m_l}<0 \text{ and } \frac{d^2}{dx^2} \big \{ \xi''(x)^{-\frac{1}{2}} \big \} \big|_{x=M_l}>0.
\end{eqnarray*}
By intermediate value theorem, there exists zeros in intervals $\big[M_{l-1},m_l \big]$ and $\big[m_l,M_l \big]$ respectively.
Since $\frac{d^2}{dx^2} \big \{ \xi''(x)^{-\frac{1}{2}} \big \} $ has at most $2(n-1)$ strictly positive roots, $\frac{d^2}{dx^2} \big \{ \xi''(x)^{-\frac{1}{2}} \big \} $ has exactly 1 root in $\big[M_{l-1},m_l \big]$ and $\big[ m_l,M_l]$ for $l=1,\cdots,n-1.$

 By \eqref{bound} and Lemma \ref{proph}, it holds that $0<q_1<M_1.$ Also for $l=2,\cdots, n-1,$ it yields that $q_{l-1}<m_l$ and $q_l<M_l.$  The interval $[q_{l-1},q_l]$ contains at most two roots of $\frac{d^2}{dx^2} \big \{ \xi''(x)^{-\frac{1}{2}} \big \} .$

 Recall for $l=1,\cdots,n$, the critical points $x \in [q_{l-1},q_l] $ of $g(x)$ satisfy that $1+z_n=\big(\sA^{n+1}_l(x)\big)^{(-1)^{n-l}}$. We then replace the system of equations in Theorem \ref{nrsb} by 
 \begin{eqnarray*}\left\{
\begin{array}{lcl}
h^n_1(\bar{q}, \sA^{n+1}_l(x))=0,\\
\qquad \cdots \cdots\\
h^n_n(\bar{q}, \sA^{n+1}_l(x))=0,
\end{array} \right . 
\end{eqnarray*}

Then in order to show the condition (3) in Theorem \ref{CriterionKRSB}, it is enough for us to show that there is exactly one $x \in (q_{l-1},q_l)$ correspondingly with $q_1, \cdots, q_{n-1}$ satisfying the system of equations above.

Define $c_l(x):=h^n_l\big(\bar{q}, (\sA^{n+1}_l(x))^{(-1)^{k-l}}\big)$. By computation, it yields that
\begin{eqnarray*}
&&\text{sgn} \Big( \frac{d}{d x} c_l(x) \Big)\\
&=&\text{sgn} \Big( (q_l-q_{l-1})[\xi'(q_l)-\xi'(x)][\xi'(x)-\xi'(q_{l-1})] -\xi''(x)[\xi'(q_l)-\xi'(q_{l-1})](q_l-x)(x-q_{l-1}) \Big)\\
&=:& \text{sgn} \Big( t_l(x) \Big),
\end{eqnarray*}
which implies that
\begin{eqnarray*}
t''_l(q_l)&=&\xi'''(q_l)[\xi'(q_l)-\xi'(q_{l-1})](q_l-q_{l-1})+2\xi''(q_l) [\xi'(q_l)-\xi'(q_{l-1})-\xi''(q_l)(q_l-q_{l-1})] \\
&=& \frac{\partial}{\partial q_l} \big( h^F(q_{l-1},q_l) \big) \\
&=& \frac{\partial}{\partial q_l} \big(  h^n_{l}(\bar{q}, \dF^n_{l}) \big),
\end{eqnarray*}
and
\begin{eqnarray*}
t_l''(q_{l-1})&=&-\xi'''(q_{l-1})[\xi'(q_l)-\xi'(q_{l-1})](q_l-q_{l-1})+2\xi''(q_{l-1})[\xi'(q_l)-\xi'(q_{l-1})-\xi''(q_{l-1})(q_l-q_{l-1})] \\
&=& -\frac{\partial}{\partial q_{l-1}} \big( h^F(q_l,q_{l-1}) \big)\\
&=& \frac{\partial}{\partial q_{l-1}}  \big(  h^n_{l}(\bar{q}, \dF^n_{l-1}) \big).
\end{eqnarray*}
Since the local extrema of $h^n_l(\bar{q}, \dF^n_{l-1}) $ and $ h^n_l(\bar{q}, \dF^n_l)$ are determined by $\frac{d^2}{dx^2} \big \{ \xi''(x)^{-\frac{1}{2}} \big \} $ and  the interval $[q_{l-1},q_l]$ contains at most two roots of $\frac{d^2}{dx^2} \big \{ \xi''(x)^{-\frac{1}{2}} \big \} $. Then $h^n_l(q_{l-1},x,\dF^n_l)$ has at most two critical points in $[q_{l-1},q_{l}]$, which implies that there are at most two zeros  of $t_l(x)$ in $[q_{l-1},q_l].$ If $t_l(x)$ has two zeros $a_l<b_l$ in $[q_{l-1},q_l]$, then it holds that $t'(a_l)<0$ and $t'(b_l)>0.$

Notice that for $n-l$ even, it holds that $\sA^{n+1}_l(q_{l-1})=(\dF^n_{l-1})^{-1}$ and  $\sA^{n+1}_l(q_l)=\dF^n_{l}.$ Also for $n-l$ odd, it holds that $\sA^{n+1}_l(q_{l-1})=\dF^n_{l-1}$ and  $\sA^{n+1}_l(q_l)=(\dF^n_{l})^{-1}.$
 By \eqref{bound} and Lemma \ref{proph}, we obtain that $c_l(q_{l-1}) \geq 0$ and $c_l(q_{l}) \leq 0.$

Therefore combining with the fact that there are at most two zeros  of $t_l(x)$ in $[q_{l-1},q_l],$ $c_l(x)$ has at most one zero in $[q_{l-1},q_l].$

\end{proof}

When $k=n,$ Theorem \ref{CriterionKFRSB} can also be simplified as follows:
\begin{theorem}\label{nfrsb}
The Parisi measure $\nu_P$ for $\xi$ is $n-$FRSB if and only if for $j=1,\cdots,t,$ $q^j_{w_{j-1}},q^j_{w_{j-1}+1},\cdots, q^j_{w_j}$ satisfy the following conditions
\begin{enumerate}
\item the system of equations
\begin{eqnarray*}
\left \{ \begin{array}{lcl}
h^j_l \Big( \bar{q}^j, *_j\Big) =0, \text{ for }l=w_{j-1}+1,\cdots,w_j,\\
\nF^{j}_{w_{j-1}}=\big(\nF^{j}_{w_j}\big)^{(-1)^{s_j}}, \text{ if } 1<j<t.
\end{array} \right.
\end{eqnarray*}
where $*_j=\nF^j_{w_j}$ if $j=1,\cdots,t-1$ and $*_t=(\nF^t_{w_{t-1}})^{(-1)^{s_t}}.$
\item $ \nY^j \leq *_j \leq (\nZ^j)^{-1} $
\end{enumerate}

\end{theorem}

\begin{proof}
By Theorem \ref{CriterionKFRSB}, it suffices for us to show that condition (3) and (4) can be deduced from conditions (1) and (2).

As for condition (3), the same argument apply  just as the proof of Theorem \ref{nrsb} based on conditions (1) and (2).

It then remains for us to show that $\frac{d^2}{dx^2} \big \{ \xi''(x)^{-\frac{1}{2}} \big \}<0$ for $x \in \big[q^j_{w_j},q^{j+1}_{w_j} \big]$ and $j=1,\cdots,t-1.$

Based on the argument in the proof of Theorem \ref{nrsb}, $\frac{d^2}{dx^2} \big \{ \xi''(x)^{-\frac{1}{2}} \big \} $ has exactly 1 root in $\big[ m_{w_j},M_{w_j}]$ for $j=1,\cdots t.$ Denote the first local maximum of $-h^F\big(q^{j+1}_{w_j+1},x\big)$  larger than $q^{j+1}_{w_j}$ by $N_{w_j}$. It holds that $\frac{d^2}{dx^2} \big \{ \xi''(x)^{-\frac{1}{2}} \big \} \big|_{x=N_{w_j}}<0.$
It then yields that $m_{w_j}<q^j_{w_j}< q^{j+1}_{w_j}<N_{w_j}<M_{w_j},$ which implies that $\frac{d^2}{dx^2} \big \{ \xi''(x)^{-\frac{1}{2}} \big \} <0,$ for $x \in [q^j_{w_j},q^{j+1}_{w_j}].$

\end{proof}

Lastly, we prove Corollary \ref{atmostnrsb} as follows:

\begin{proof}[Proof of Corollary \ref{atmostnrsb}]
Based on the proof of Theorem \ref{nrsb} and \ref{nfrsb}, $\frac{d^2}{dx^2} \big \{ \xi''(x)^{-\frac{1}{2}} \big \}$ must have at least $2(k-1)$ strictly positive roots so that $\nu_P$ is $k-RSB.$ Since for $\xi(x)=\sum^{n-1}_{l=1} \lambda_l \cdot x^{p_l} +\Big( 1- \sum^{n-1}_{l=1} \lambda_l \Big) \cdot x^{p_n},$ $\frac{d^2}{dx^2} \big \{ \xi''(x)^{-\frac{1}{2}} \big \}$ has at most $2(n-1)$ strictly positive roots, we conclude that $k \leq n.$

\end{proof}

\section{Proof of Theorem \ref{kRSB},\ref{kFRSB} and Corollary \ref{nRSB},\ref{nFRSB}}\label{Secproofmain}

Recall the notations in Section \ref{main}, for $j=1,\cdots,t$ and any $\bar{q}^j=(q^j_{w_{j-1}},q^j_{w_{j-1}+1},\cdots,q^j_{w_j}) \in \mathcal{H}^{s_j},$ it holds that
\begin{eqnarray*} h^j_l(\bar{q}^j, (\mathcal{Z}^j)^{-1})= 
\left \{ \begin{array}{lcl}
h\Big(q^j_{l-1},q^j_l, \big({\nF^j_l} \cdot \mathcal{Z}^j \big)^{-1}\cdot r_2(q^j_{l-1},q^j_l) \Big)\geq 0, \text{ for } w_j-l \text{ even,}\\
h\Big(q^j_{l-1},q^j_l, \frac{\mathcal{Z}^j}{\nF^j_l}\cdot r_2(q^j_{l-1},q^j_l) \Big) \leq 0 , \text{ for } w_j-l \text{ odd},
\end{array} \right.
\end{eqnarray*}
and
\begin{eqnarray*} h^j_l(\bar{q}^j, \mathcal{Y}^j)=
\left \{ \begin{array}{lcl}
h \Big(q^j_{l-1},q^j_l, \frac{\mathcal{Y}^j}{\nF^j_l} \cdot r_2(q^j_{l-1},q^j_l)\Big)  \leq 0 , \text{ for } w_j-l \text{ even,}\\
h \Big(q^j_{l-1},q^j_l, \big(\mathcal{Y}^j \cdot {\nF^j_l} \big)^{-1} \cdot r_2(q^j_{l-1},q^j_l)\Big)  \geq 0  , \text{ for } w_j-l \text{ odd}.
\end{array} \right.
\end{eqnarray*}

In particular for $\bar{q}=(0,q_1,\cdots,q_{k-1},1) \in \mathcal{H}^k,$ it holds that for $l=1,\cdots,k,$
\begin{eqnarray*} h^k_l(\bar{q}, (\mathcal{Z}^k)^{-1})= 
\left \{ \begin{array}{lcl}
h\Big(q_{l-1},q_l, \frac{\dF^k_{l-1}(\bar{q})}{\mathcal{Z}^k(\bar{q})} \cdot r_2(q_l,q_{l-1})^{-1} \Big) \geq 0, \text{ for } k-l \text{ even,}\\
h\Big(q_{l-1},q_l, \frac{\mathcal{Z}^k(\bar{q})}{\nF^k_l(\bar{q})}\cdot r_2(q_{l-1},q_l) \Big) \leq 0 , \text{ for } k-l \text{ odd},
\end{array} \right.
\end{eqnarray*}
and
\begin{eqnarray*} h^k_l(\bar{q}, \mathcal{Y}^k)=
\left \{ \begin{array}{lcl}
h \Big(q_{l-1},q_l, \frac{\mathcal{Y}^k(\bar{q})}{\dF^k_l(\bar{q})} \cdot r_2(q_{l-1},q_l)\Big)  \leq 0 , \text{ for } k-l \text{ even,}\\
h\Big(q_{l-1},q_l, \frac{\dF^k_{l-1}(\bar{q})}{\mathcal{Y}^k(\bar{q})} \cdot r_2(q_l,q_{l-1})^{-1} \Big) \geq 0 , \text{ for } k-l \text{ odd}.
\end{array} \right.
\end{eqnarray*}

 We then prove Theorem \ref{kRSB} as follows:
 \begin{proof}[Proof of Theorem \ref{kRSB}]
 It is equivalent for us to show that conditions (1),(2) and $(3)$ in Theorem \ref{CriterionKRSB} are equivalent to conditions (1),(2) in Theorem \ref{kRSB}.
 
 We first prove that conditions (1),(2)  in Theorem \ref{CriterionKRSB} is equivalent to condition (1) in Theorem \ref{kRSB}.
 
 Assume that conditions (1),(2) in Theorem \ref{CriterionKRSB} hold.
 Since for $l=1,\cdots, k$, it holds that, if $k-l$ even,
 \begin{eqnarray}
h\Big(q_{l-1},q_{l}, (1+z_k) \dF^k_{l-1} \cdot r_2(q_{l},q_{l-1})^{-1}  \Big) =h\Big(q_{l-1},q_{l}, \big(\frac{1+z_k}{\dF^k_l} \big) \cdot r_2(q_{l-1},q_l)  \Big) =0,
\end{eqnarray}
and if $k-l$ odd,
\begin{eqnarray}\label{1+zodd}
h\Big(q_{l-1},q_{l}, \big(\dF^k_{l}(1+z_k)\big) ^{-1}  \cdot  r_2(q_{l-1},q_{l})\Big)=h\Big(q_{l-1},q_{l}, \big(\frac{\dF^k_{l-1}} {1+z_k}\big)  \cdot  r_2(q_{l},q_{l-1})^{-1} \Big)=0.
\end{eqnarray}
which is equivalent to $h^k_l(\bar{q},1+z_k)=0$ for $k=1,\cdots,l.$
 
Also by the assumption $\mathcal{Y}^k  \leq 1+z_k \leq \big(\mathcal{Z}^k\big)^{-1}$ and Lemma \ref{proph},   it then yields that $\bar{q}=(0, q_1,q_2, \cdots, q_{k-1},1) $ satisfy Condition $\varkappa(k)$ and then $\bar{q} \in \mathcal{H}^k ,$ which implies condition (1) in Theorem \ref{kRSB}.

 Now we assume condition (1) in Theorem \ref{kRSB} holds. It is then implied, by intermediate value theorem, there exists $0 < q_1<q_2<\cdots<q_{k-1}< 1 $ and $1+z_k \in [\mathcal{Y}^k,\big(\mathcal{Z}^k\big)^{-1}]$ satisfying that \eqref{1+zeven} and $\eqref{1+zodd}$ for $l=1,\cdots, k$.

 Nest we assume that conditions (1),(2) hold in Theorem \ref{kRSB}  and prove condition $(3)$ in Theorem \ref{CriterionKRSB}.


For $l=1,\cdots,k$, assume that $x^*_l$ is a critical point of $g(x)$ in $[q_{l-1},q_l]$ and set $\bar{x}_l=(0,q_1,\cdots,q_{l-1},x^*_l,q_l,\cdots,q_{l-1},1)$. Without loss of generality, we assume $x^*_l$ is a local minimum of $g(x)$, which implies that
$\sA^{k+1}_l(x^*_l)>\dF^{k+1}_l(x^*_l)$ for $k-l$ odd 
and $\dF^{k+1}_{l}(x^*_l)>\sA^{k+1}_l(x^*_l)$ for $k-l$ even.

We notice that for $l'=l+1,\cdots,k,$ $$h^k_{l'}(\bar{q},1+z_k)=h^{k+1}_{l'+1}(\bar{x}_l,1+z_k)=0,$$
and for $l'=1,\cdots,l-1,$
$$h^k_{l'}(\bar{q},1+z_k)=h^{k+1}_{l'}(\bar{x}_l,1+z_k)=0 \text{ if } 1+z_k=\big(\sA^{k+1}(x^*_l)\big)^{(-1)^{k-l}}.$$

Also note that for any $x \in [q_{l-1},q_l],$
 \begin{eqnarray*}
 &&h^k_l\big(\bar{q}, \sA^{k+1}_{l}(x)\big),\\
 &=&h^{k+1}_l(\bar{x}_l,\sA^{k+1}_{l}(x))+h^{k+1}_{l+1}(\bar{x}_l,\sA^{k+1}_{l}(x)).
 \end{eqnarray*}
Recall that $h^k_l\big(\bar{q}, \sA^{k+1}_{l}(x^*_l)\big)=0.$ Then if it holds that $h^{k+1}_{l+1}\big(\bar{x}_l, \sA^{k+1}_{l}(x^*_l)\big)<0$, then $h^{k+1}_{l} \big(\bar{x}_l,\sA^{k+1}_{l}(x^*_l)\big)>0. $ 
Therefore $\bar{x}_l=(0,q_1,\cdots,q_{l-1},x^*_l,q_l,q_{k-1},1) \in \mathcal{H}^{k+1}$, which contradicts with condition (2) in Theorem \ref{kRSB}. Then condition (3) in Theorem \ref{CriterionKRSB} holds.

 Next, we assume condition (3) in Theorem \ref{CriterionKRSB} holds. We prove condition (2) in Theorem \ref{kRSB} by contradiction and assume $\mathcal{H}^{k+1} \neq \emptyset.$ If $\mathcal{H}^{k+2}= \emptyset$, then by previous reasoning, the Parisi measure $\nu_P$ is then $k+1-$RSB, which leads to a contradiction. If $\mathcal{H}^{k+2}\neq \emptyset$ and $\mathcal{H}^{k+3}= \emptyset,$ the same argument applies. By Theorem \ref{nrsb}, it always holds that $\mathcal{H}^{n+1}= \emptyset,$ which finishes the proof.


 

 \end{proof}
 
 Now we turn to the proof of Theorem \ref{kFRSB}.
 \begin{proof}[Proof of Theorem \ref{kFRSB}]
  It is equivalent for us to show that conditions (1-4) in Theorem \ref{CriterionKFRSB} are equivalent to conditions (1-2) in Theorem \ref{kFRSB}.

We first assume conditions (1-2) in Theorem \ref{kFRSB} hold.

Since $\mathcal{H}^{s_1} \tilde{<} \mathcal{H}^{s_2} \tilde{<} \cdots \tilde{<} \mathcal{H}^{s_t} $, then for $j=1,\cdots,t-1$ and $\bar{q}^j=(q_{w_{j-1}},q_{w_{j-1}+1},\cdots,q_{w_j}) \in \mathcal{H}^{s_j}$ satisfying that  $q_{w_j}=\max_{(r_0,r_1\cdots,r_{s_j})\in \mathcal{H}^k_{s_j}} r_{s_j}$, it holds that $\mathcal{Y}^{s_j}=\nF^j_{w_j}$ and $h^j_l(\bar{q}^j,\nF^j_{w_j})=0$, $l=w_{j-1}+1,\cdots,w_j.$

Also for $j=2,\cdots,t$ and $\bar{q}^j=(q_{w_{j-1}},q_{w_{j-1}+1},\cdots,q_{w_j}) \in \mathcal{H}^{s_j}$ satisfying that $q_{w_{j-1}}=\min_{(r_0,r_1\cdots,r_{s_j})\in \mathcal{H}^k_{s_j}} r_{0}$, it holds that $\nF^j_{w_{j-1}}=\mathcal{Y}^{s_j}$ if $s_j$ is even and $\nF^j_{w_{j-1}}=\mathcal{Z}^{s_j}$ if $s_j$ is odd. Moreover, $h^j_l \big(\bar{q}^j,(\nF^j_{w_{j-1}})^{(-1)^{s_j}} \big)=0$, $l=w_{j-1}+1,\cdots,w_j.$

Furthermore for $j=1,\cdots, t$ and $\bar{q}^j=(q_{w_{j-1}},q_{w_{j-1}+1},\cdots, q_{w_j} )\in \mathcal{H}^{s_j} $,
\begin{eqnarray*}
 \text{sgn} \Big( \frac{\partial}{\partial q_{w_j-1}} h^j_{w_j}(\bar{q}^j, F^j_{w_j})  \Big)
 = \text{sgn} \big( 1-\nF^j_{w_j} \cdot \nF^j_{w_j-1} \big) \geq0, 
\end{eqnarray*}
and
\begin{eqnarray*}
 \text{sgn} \Big( \frac{\partial}{\partial q_{w_{j-1}+1}} h^j_{w_{j-1}+1} \big(\bar{q}^j ,(F^j_{w_{j-1}})^{(-1)^{s_j}} \big)  \Big)
 = 
\text{sgn} \big( \nF^j_{w_{j-1}+1} \cdot \nF^j_{w_{j-1}} -1 \big) \leq0. 
\end{eqnarray*}
Then when $q_{w_{j-1}+1}$ is increasing, the zero of $h^j_{w_{j-1}+1} \big(\bar{q}^j,(F^j_{w_{j-1}})^{(-1)^{s_j}} \big) $ is also increasing. Also when $q_{w_j-1}$ is increasing, the zero of $h^j_{w_j}(\bar{q}^j,F^j_{w_j}) $ is decreasing.

Therefore by condition (1) in Theorem \ref{kFRSB}, there exists $\bar{q}^j =(q^j_{w_{j-1}},\cdots,q^j_{w_j})\in \mathcal{H}^{s_j}$ such that $q^{j}_{w_j}= \max_{(r_0,r_1\cdots,r_{s_j})\in \mathcal{H}^{s_j}} r_{s_j}$ and $q^{j}_{w_{j-1}}=\min_{(r_0,r_1\cdots,r_{s_j})\in \mathcal{H}^{s_j}} r_{0}$ with $q^{j}_{w_j}<q^{j+1}_{w_j}$, for $j=1,\cdots,t-1.$ Moreover $\bar{q}^j$ also satisfy that $\mathcal{Y}^{j} \cdot \mathcal{Z}^{j} <1$  and the system of equations \eqref{conditionkfrsb}, which implies condition (1-2) in Theorem \ref{CriterionKFRSB}.

Since there exists no $j=1,\cdots,t$ such that the relation $\mathcal{H}^{s_1} \cdots \tilde{<} \mathcal{H}^{s_{j-1}} \tilde{<}  \mathcal{H}^{1} \tilde{<}  \mathcal{H}^{s_{j}} \tilde{<} \cdots \tilde{<} \mathcal{H}^{s_t} $ hold. Then for any $(x,y)$, it holds that $h^F(x,y) \geq 0$ and $h^F(y,x) \geq 0.$
which implies that $\frac{d^2}{dx^2} \big \{ \xi''(x)^{-\frac{1}{2}} \big \} \leq 0$, for $x \in [q^j_{w_j},q^{j+1}_{w_j}].$

Assume that $x^*_l$ is a critical point of $g(x)$ in $[q^j_{l-1},q^j_l]$. Without loss of generality, we assume $x^*_l$ is a local minimum of $g(x)$, which implies that
$\nA_l^j(x^*_l)\cdot \nF^j_l \leq 1.$

Now we define $\tilde{h}^j_l(\bar{x}^j_l, *_j)=g(q^j_{l-1})$ for $s=s_j+1$ and $l = w_{j-1}+1, \cdots w_j+1$.

We notice that for $l'=l+1,\cdots,w_j,$ $$h^j_{l'}(\bar{q}^j,*_j)=\tilde{h}^j_{l'+1}(\bar{x}^j_l,*_j)=0,$$
and for $l'=w_{j-1}+1,\cdots,l-1,$
$$h^j_{l'}(\bar{q}^j,*_j)=\tilde{h}^{j}_{l'}(\bar{x}_l,*_j)=0 \text{ if } *_j=\big(\sA^{k+1}(x^*_l)\big)^{(-1)^{w_j-l}}.$$

Also note that for any $x \in [q_{l-1},q_l],$
 \begin{eqnarray*}
 &&h^j_l\big(\bar{q}, \nA^{j}_{l}(x)\big),\\
 &=&\tilde{h}^{j}_l(\bar{x}_l,\nA^{j}_{l}(x))+\tilde{h}^{j}_{l+1}(\bar{x}_l,\nA^{j}_{l}(x)).
 \end{eqnarray*}
Recall that $h^j_l\big(\bar{q}, \nA^{j}_{l}(x^*_l)\big)=0.$ Then if it holds that $h^{j}_{l+1}\big(\bar{x}_l, \nA^{j}_{l}(x^*_l)\big)<0$, then $h^{k+1}_{l} \big(\bar{x}_l,\sA^{k+1}_{l}(x^*_l)\big)>0. $ 
Therefore $\bar{x}^j_l \in \mathcal{H}^{s_j+1}$ and $\mathcal{H}^{s_1} \cdots \tilde{<} \mathcal{H}^{s_{j-1}} \tilde{<}  \mathcal{H}^{s_{j}+1} \tilde{<}  \mathcal{H}^{s_{j+1}} \tilde{<} \cdots \tilde{<} \mathcal{H}^{s_t} $, which contradicts with condition (2) in Theorem \ref{kFRSB}. Then condition (3) in Theorem \ref{CriterionKFRSB} holds.

Now assume conditions (1-4) in Theorem \ref{CriterionKFRSB} holds.
Since for $j=1,\cdots,t+1,$ $\bar{q}^j=(q^j_{w_{j-1}},q^j_{w_{j-1}+1},\cdots, q^j_{w_j})$ satisfy that  $h^j_l(\bar{q}^j, \nF^j_{w_{j}})=0,$ where $l=w_{j-1}+1,\cdots,w_j$ and $\nY^j \cdot \nZ^j \leq 1$, it yields that $\bar{q}^j \in \mathcal{H}_{s_j} $.

For any $(x_l,x_{l+1}) \in \mathcal{H}^{1}$, it holds that $h^F(x_l,x_{l+1}) \leq 0$ and $h^F(x_{l+1},x_l) \leq 0.$
Since $\frac{d^2}{dx^2} \big \{ \xi''(x)^{-\frac{1}{2}} \big \} \leq 0$, for $x \in [q^j_{w_j},q^{j+1}_{w_j}],$ it yields that $h^F(x_l,x_{l+1})>0$ and $h^F(x_{l+1},x_l)<0$ for any $q^j_{w_j}<x_l<x_{l+1}<q^{j+1}_{w_j}.$ Thus there exists no $\mathcal{H}^{1} $ such that $\mathcal{H}^{s_j} \tilde{<} \mathcal{H}^{1} \tilde{<} \mathcal{H}^{s_{j+1}}$.

Next we assume that there exists $1 \leq j \leq t$ such that the relation $\mathcal{H}^{s_1} \cdots \tilde{<} \mathcal{H}^{s_{j-1}} \tilde{<}  \mathcal{H}^{s_{j}+1} \tilde{<}  \mathcal{H}^{s_{j+1}} \tilde{<} \cdots \tilde{<} \mathcal{H}^{s_t} $ holds. By Theorem \ref{nfrsb}, without loss of generality, we assume that the relation $\mathcal{H}^{s_1} \cdots \tilde{<} \mathcal{H}^{s_{j-1}} \tilde{<}  \mathcal{H}^{s_{j}+2} \tilde{<}  \mathcal{H}^{s_{j+1}} \tilde{<} \cdots \tilde{<} \mathcal{H}^{s_t} $ doesn't hold.  Then by previous reasoning, the Parisi measure $\nu_P$ is $(k+1)$-FRSB, which leads to a contradiction.


 \end{proof}

Finally we prove Corollary \ref{nRSB} and \ref{nFRSB} as follows:
\begin{proof}[Proof of Corollary \ref{nRSB} and \ref{nFRSB}]

Since by Theorem \ref{nrsb}, condition (2) in Theorem \ref{kRSB} and equivalently condition (3) in Theorem \ref{CriterionKRSB} are automatically satisfied when $k=n,$ which then implies Corollary \ref{nRSB}.

Similarly by Theorem \ref{nfrsb}, condition (2),(3) in Theorem \ref{kFRSB} and equivalently condition (3)(4) in Theorem \ref{CriterionKFRSB} are automatically satisfied when $k=n,$ which then implies Corollary \ref{nRSB}.

\end{proof}

\bibliographystyle{abbrv}
\bibliography{biblio3}

\end{document}